\documentclass[12pt]{article}
\usepackage{a4, amssymb, latexsym, amsmath, exscale, amscd}
\usepackage[all]{xy}

\newenvironment{evlist}[2]{
\begin{list}{}{
\setlength{\topsep}{0.5ex plus0.2ex minus0.1ex} 
\setlength{\leftmargin}{#1}
\setlength{\itemsep}{#2 plus0.2ex}
\setlength{\parsep}{0ex plus0.2ex} }}
{\end{list}}

\newcommand{\mtimes}{\times}
\newcommand{\btimes}{\diamond}
\newcommand{\Fam}[1]{\mathrm{#1}}
\newcommand{\Self}[1]{\mathrm{T}_{#1}}
\newcommand{\lsm}[1]{\langle{#1}\rangle}
\newcommand{\Nat}{\mathbb{N}}
\newcommand{\Int}{\mathbb{Z}}
\newcommand{\id}{\mathrm{id}}
\newcommand{\Hom}{\mathrm{Hom}}

\newcommand{\definition}[1]{\textit{#1}}

\newcommand{\proof}{{\textit{Proof}\enspace}}

\newcommand{\eop}{\ \vbox{\hrule
                       \hbox{\vrule
                             \hskip 6pt
                             \vrule height 6pt width 0pt
                             \vrule}%
                       \hrule}%
                     \vspace{\medskipamount}
                }

\newsavebox{\ttt}
\sbox{\ttt}{}
\newcommand{\startsection}[1]
    {\section[#1]{#1}
    \sbox{\ttt}{\thesection\ \ \textsc{#1}}
    \thispagestyle{plain}
}

\newtheorem{lemma}{Lemma}[section]
\newtheorem{proposition}{Proposition}[section]
\newtheorem{theorem}{Theorem}[section]

\newcommand\mythicklines{}
\newcommand\sput[3]{\put(#1,#2){$\scriptstyle{#3}$}}

\newlength{\graphicthick}
\newlength{\graphicmid}
\newlength{\graphicthin}

\begin{document}

\begin{titlepage}

\begin{center}
\phantom{q}
\vspace{60pt}

{\huge Minimal counting systems}

{\LARGE and }

\smallskip
{\huge commutative monoids}

\bigskip\bigskip\bigskip
{\LARGE Chris Preston}

\bigskip
%{\large February 2009}
\end{center}

\bigskip
\bigskip
\bigskip
\bigskip
\begin{quote}
These notes present an approach to obtaining the basic operations of addition 
and multiplication on the natural numbers in terms of elementary results
about commutative monoids.
\end{quote}

\end{titlepage}

\thispagestyle{empty}

\addtocontents{toc}{\vskip 20pt}
\addtolength{\parskip}{-5pt}
\tableofcontents
\addtolength{\parskip}{5pt}

\startsection{Introduction}

\label{intro}

These notes arose from a  first semester course given 
to students studying to be primary school mathematics teachers.
One aim of the course was to explain how the operations of addition and multiplication
can be introduced within the framework of the Peano axioms.
Perhaps the most natural approach  is in terms of
elementary operations on finite sets (disjoint union and cartesian product)
and in Preston \cite{preston} we presented the mathematics behind
one particular version of this approach. The material here and in \cite{preston}
are aimed at those who have to teach such courses and  certainly not meant
for the students; they could also be used as the basis for a `proper' mathematics seminar.

The `natural' interpretation of $m + n$ is that it is the number of elements in the disjoint union
of two finite sets, one containing $m$ elements and the other $n$.
In particular, $m$ and $n$ have equal roles here, and $+$ is thought of as an operation
with two arguments. But there is also the `operator' interpretation, in which  $m + n$ is regarded
as the result of applying the operation $+\,n$ to $m$. 
Here $m$ and $n$ have completely different roles, with $m$ appearing as the argument
of a mapping $+\,n : \Nat \to \Nat$.

It is the mappings $+\,n$, $n \in \Nat$, which form the basis of the present approach.
These mappings lie in $\Self{\Nat}$, the set of all mappings of $\Nat$ into itself, which is
a monoid with functional composition $\circ$ as monoid operation
and $\id_\Nat$ as unit element. 
Now the set $M = \{ +\,n : n \in \Nat \}$ is a commutative submonoid of $\Self{\Nat}$
and the mapping $+\,(m + n)$ in $M$ is just the functional composition of the mappings $+\,m$ and $+\,n$.
This means that the monoids $(M,\circ,\id_\Nat)$ and $(\Nat,+,0)$ are isomorphic,
and in fact the evaluation mapping $\Phi_0 : M \to \Nat$ with $\Phi_0(u) = u(0)$ for all
$u \in M$ is an isomorphism.

These observations require, of course, that we already know what the addition is. However, the above 
procedure can be reversed and, starting with the data given by the Peano axioms, it can be used to define 
$+$. This is because the mapping $+\,1$ is nothing but the successor operation
$\textsf{s} : \Nat \to \Nat$  (with $\textsf{s}(0) = 1$, $\textsf{s}(1) = 2$ and so on)
and the submonoid $M$ is generated by the single mapping $+\,1$. 
The addition can therefore by obtained as follows: First define $M$ to be the submonoid of
$\Self{\Nat}$ generated by $\textsf{s}$; such a submonoid is always commutative.
Then show that the evaluation mapping $\Phi_0 : M \to \Nat$ is a bijection.
Finally, use the bijection $\Phi_0$ to transfer the monoid operation $\circ$ on $M$
to an operation $+$ on $\Nat$, which is then the required addition.

What does this have to do with the Peano axioms? For the construction to work we require only that
$\Phi_0$ be a bijection, which turns out to be  equivalent to the principle
of mathematical induction being valid. This means that an addition can be defined without 
the other Peano axioms holding. For example, the addition on $\Int_n$
can be obtained directly in this fashion.

The method also works when the single successor mapping $\textsf{s}$ is replaced by a family of commuting
mappings. For example, the addition on the integers $\Int$ can be obtained by using
the successor operation together with its inverse.

Using this approach the multiplication $\times$ will be obtained from a simple result 
for commutative monoids which guarantees the existence of biadditive mappings 
specified on a set of generators of the monoid.

Another topic which fits well into the present approach concerns generalisations of the recursion theorem.
These will be dealt with as problems involving free commutative monoids.

%%% Local Variables: 
%%% mode: latex
%%% TeX-master: "sums"
%%% End: 

\startsection{Sums and products of natural numbers}

\label{spnat}

In this section we do the following:

\begin{evlist}{18pt}{6pt}
\item[--]
Start with a reminder of how the operations of addition and multiplication for the natural numbers are
usually presented within the framework of the Peano axioms.

\item[--]
Point out that these operations can be defined without having 
all of the axioms; all that is needed is that the principle of mathematical induction should hold.

\item[--]
Explain how the operations can be  obtained with the help of some simple
constructions involving commutative monoids. 
\end{evlist}

Some of the details of the proofs are left to the following sections and are carried out in a more
general set-up. This arises by replacing the single successor mapping occurring in the Peano axioms
with a family of commuting mappings.

Consider the set of natural numbers $\Nat = \{0,1,\ldots\}$ together with the successor operation
$\textsf{s} : \Nat \to \Nat$  (with $\textsf{s}(0) = 1$, $\textsf{s}(1) = 2$ and so on). 
Then $(\Nat,\textsf{s},0)$ is an example of a
triple $(X,f,x_0)$ consisting of a set $X$, a mapping $f : X \to X$ and an element $x_0 \in X$, and we call
any such triple a \definition{counting system}. 
The question of how the special
counting system $(\Nat,\textsf{s},0)$ should be characterised 
within the class of all counting systems was answered by Dedekind in his book 
\textit{Was sind und was sollen die Zahlen?}\ \cite{dedekind}.
The requirements introduced by Dedekind to characterise $(\Nat,\textsf{s},0)$ 
are those which are now usually referred to as the
Peano axioms.

There are three axioms.
One is the \definition{principle of mathematical induction}, which says that the only 
$f$-invariant subset of $X$ containing $x_0$ is $X$ itself (and where a subset $Y \subset X$ is 
\definition{$f$-invariant} if $f(Y) \subset Y$). The other two axioms
require that the mapping $f$ be injective and that $f(x) \ne x_0$ for all $x \in X$.

A counting system for which the principle of mathematical induction holds will be called \definition{minimal}
and one obeying all the Peano axioms will be called a \definition{Dedekind system}. 
If $(X,f,x_0)$ is a Dedekind system then, as most mathematics students used to learn at some time
during their studies, there exists a unique binary operation $+$ on $X$ such that
\begin{evlist}{36pt}{6pt}
\item[$\mathrm{(+_0)}$]
$\ x_0 + x = x$ for all $x \in X$,
\item[$\mathrm{(+_1)}$]
$\ f(x_1) + x_2 = f(x_1 + x_2)$ for all $x_1,\,x_2 \in X$.
\end{evlist}
Moreover, the operation $+$ is both associative and commutative. 

\medskip
\newpage

There is also a unique binary operation $\mtimes$ on $X$ such that
\begin{evlist}{36pt}{6pt}
\item[$\mathrm{(\mtimes_0)}$]  
$\ x_0 \mtimes x = x_0$ for all $x \in X$,
\item[$\mathrm{(\mtimes_1)}$]  
$\ f(x_1) \mtimes x_2 = x_2 + (x_1 \mtimes x_2)$ for all $x_1,\,x_2 \in X$.
\end{evlist}
Again, the operation $\mtimes$ is both associative and commutative and, 
in addition, the distributative law holds for $+$ and $\mtimes$, meaning
that
\[ x_1 \mtimes (x_2 + x_3) = (x_1 \mtimes x_2) + (x_1 \mtimes x_3)\] 
for all $x_1,\,x_2,\,x_3 \in X$. It follows immediately from these properties that
$f(x_0)$ is a unit element for $\mtimes$, i.e., $f(x_0) \mtimes x = x$ for all $x \in X$.

The first proof of these facts can be found in 
\cite{dedekind}, where Dedekind
obtains the operations with the help of the recursion theorem. This theorem states that 
a Dedekind system $(X,f,x_0)$ is \definition{initial}, meaning
that for each counting system $(Y,g,y_0)$ there exists a unique mapping $\pi : X \to Y$ with 
$\pi(x_0) = y_0$ such that $\pi \circ f = g \circ \pi$. 
Being initial is probably the most useful property of a Dedekind system, since 
it is the basis for making inductive (or recursive) definitions. 
Now the recursion theorem does not hold in general for minimal counting systems,
which might suggest that the full force of the Peano axioms
is required to obtain $+$ and $\mtimes$.
However, this is not the case. As indicated at the start of the section, operations 
satisfying all the properties stated above exist as soon as 
the principle of mathematical induction holds; the other two Peano axioms are not needed.
Moreover, there are analogous results when the single mapping $f$ is replaced by a
family of commuting mappings, and among other things
it is these results we are going to present in these notes.

If $+$ is the operation described above then $(X,+,x_0)$ is a commutative monoid, 
and together with counting systems these are the main algebraic structures to be employed here. 
Let us  thus review some pertinent definitions. A \definition{monoid} 
is a triple $(M,\bullet,e)$ consisting of a set, an associative
binary operation $\bullet$ on $M$ (and so $(a_1 \bullet a_2) \bullet a_3 = a_1 \bullet (a_2 \bullet a_3)$ 
for all $a_1,\,a_2,\,a_3 \in M$) and a unit element $e \in M$ (meaning that $a \bullet e = e \bullet a = a$
for all $a \in M$). In fact $e$ is \definition{the} unit element, since
if also $a \bullet e' = e' \bullet a = a$ for all $a \in M$ then $e' = e$.
The monoid is \definition{commutative} if
$a_1 \bullet a_2 = a_2 \bullet a_1$ for all $a_1,\,a_2 \in M$, and a
subset $M'$ of $M$ is a \definition{submonoid}  if $e \in M'$ and $a_1 \bullet a_2 \in M'$
for all $a_1,\,a_2 \in M'$.

We usually just write $M$ instead of $(M,\bullet,e)$ with the operation being
denoted generically by $\bullet$ and the unit by $e$. However,
if the monoid is commutative then the operation will be denoted by $+$ and the unit, which is then called the 
\definition{zero}, by $0$ (unless it has previously been given some other denotation such as $x_0$).

If $M$ and $N$ are monoids then $\alpha : M \to N$ is a \definition{homomorphism} if $\alpha(e) = e$ 
and $\alpha(a_1 \bullet a_2) = \alpha(a_1) \bullet \alpha(a_2)$ for all $a_1,\,a_2 \in M$.
A homomorphism $\alpha : M \to N$ is an \definition{isomorphism} is there exists a homomorphism
$\alpha' : N \to M$ with $\alpha' \circ \alpha = \id_M$ and $\alpha \circ \alpha' = \id_N$, which is the
case if and only if $\alpha$ is a bijection.

Here is the main result giving the existence of the operation $+$ assuming only
that the principle of mathematical induction holds. In Proposition~\ref{prop_spnat_11} below it is shown
that $+$ is the unique operation satisfying $\mathrm{(+_0)}$ and $\mathrm{(+_1)}$.

\begin{theorem}\label{theorem_spnat_11}
Let $(X,f,x_0)$ be a minimal counting system. Then there exists a unique binary operation $+$ such that
$(X,+,x_0)$ is a commutative monoid with
\begin{evlist}{12pt}{6pt}
\item[$\mathrm{(\star)}$]
$\ f(x) = f(x_0) + x$ for all $x \in X$.
\end{evlist}
Moreover, $X$ is an abelian group if and only if $f$ is a bijection, and it 
obeys the cancellation law (meaning that $x_1 = x_2$ holds whenever $x + x_1 = x + x_2$ 
for some $x \in X$) if and only if $f$ is injective. 
\end{theorem}

\proof
This is a special case of Theorem~\ref{theorem_mcs_11}. However,
let us sketch the main part of the proof in order to give an idea of what is involved.

Denote by $\Self{X}$ the set of all mappings of $X$ into itself.
Then there is the monoid $(\Self{X},\circ,\id_X)$, where $\circ$ is functional composition and 
$\id_X$ the identity mapping.
Let $M_f$ be the least submonoid of $\Self{X}$ containing $f$; then $M_f$ is commutative 
(see Lemma~\ref{lemma_mcs_11}).
Denote by $\Phi_{x_0} : M_f \to X$ the evaluation mapping given by $\Phi_{x_0}(u) = u(x_0)$
for each $u \in M_f$. 
The crucial fact, which is established in Section~\ref{mcs}, 
is that the counting system $(X,f,x_0)$ being minimal implies the mapping $\Phi_{x_0}$ is a bijection
and so the monoid structure on $M_f$ can be transferred to $X$: There exists a unique binary operation $+$
on $X$ such that
\[ \Phi_{x_0}(u_1 \circ u_2) = \Phi_{x_0}(u_1) + \Phi_{x_0}(u_2) \] 
for all $u_1,\,u_2 \in M_f$. With this operation $X$ becomes a commutative monoid
with zero element $x_0 = \Phi_{x_0}(\id_X)$, and $\Phi_{x_0}$ is then an isomorphism of monoids.
\[ \begin{CD}
M_f \times M_f @>\circ >> M_f \\
@V\Phi_{x_0}\times \Phi_{x_0} VV              @VV\Phi_{x_0} V       \\
  X \times X @>> + > X  
\end{CD} \]
Now to see that $\mathrm{(\star)}$ holds, let $x \in X$ and $u \in M_f$;
since $\Phi_{x_0}$ is surjective there exists $v \in M_f$
with $x = \Phi_{x_0}(v) = v(x_0)$. Then
\[ u(x) = u(v(x_0)) = \Phi_{x_0}(u \circ v)  = \Phi_{x_0}(u) + \Phi_{x_0}(v)  
= \Phi_{x_0}(u) + x = u(x_0) + x \]
and $\mathrm{(\star)}$ is the special case with $u = f$. 
The details and the remaining parts of the proof are contained in the proof of Theorem~\ref{theorem_mcs_11}.
\eop

The operation $+$ in Theorem~\ref{theorem_spnat_11} will be called the
\definition{operation associated with $(X,f,x_0)$} and
the resulting monoid $X$ will  be referred to as the \definition{associated monoid}.

\begin{proposition}\label{prop_spnat_11}
The operation $+$ associated with the minimal counting system $(X,f,x_0)$ is  
the unique binary operation on $X$ 
satisfying $\mathrm{(+_0)}$ and $\mathrm{(+_1)}$.
\end{proposition}

\proof
It is clear that $\mathrm{(+_0)}$ holds, and
from $\mathrm{(\star)}$ it follows that
\[ f(x_1) + x_2 = (f(x_0) + x_1) + x_2 = f(x_0) + (x_1 + x_2) = f(x_1 + x_2)\]
for all $x_1,\,x_2 \in X$, i.e., $\mathrm{(+_1)}$ holds.
Let $+'$ be a further operation satisfying $\mathrm{(+_0)}$ and $\mathrm{(+_1)}$ and put
$E = \{ y \in X : x +' y = x + y\ \mbox{for all $x \in X$} \}$.
Then $x_0 \in E$, since
$x_0 + x  = x = x_0 +' x$ for all $x \in X$, and if $y \in E$ then
$f(y) \in E$, since
$f(y) + x = f(y + x) = f(y +' x) = f(y) +' x$ for all $x \in X$.
Thus $E$ is an $f$-invariant subset of $X$ containing $x_0$, and hence $E = X$, 
which implies that ${+'} = {+}$. \eop

The following additional property of $+$ plays an important role when
introducing the multiplication $\times$:

\begin{lemma}\label{lemma_spnat_11}
If $(X,f,x_0)$ is a minimal counting system then the single element $f(x_0)$
generates the associated monoid (meaning that the only submonoid of $X$ containing $f(x_0)$ is $X$ itself).
\end{lemma}

\proof 
Let $X'$ be any submonoid of $X$ containing $f(x_0)$. If $x \in X'$ then by $\mathrm{(\star)}$
$f(x) = f(x_0) + x \in X'$, and thus $X'$ is an $f$-invariant
subset of $X$ which contains $x_0$ (since $X'$ is a submonoid). Thus $X' = X$.
\eop 

The second operation $\times$ will be obtained via a result
involving a commutative monoid generated by a single 
element, which by Lemma~\ref{lemma_spnat_11}
can be applied to the monoid associated with a minimal counting system.
Let $M$ be any commutative monoid.
A mapping $\btimes : M \times M \to M$ is said to be \definition{biadditive} if
$a' \mapsto a \btimes a'$ and $a' \mapsto a' \btimes a$ are both endomorphism of $M$ 
for each $a \in M$.
(An \definition{endomorphism} $\alpha$ of $M$ is a homomorphism $\alpha : M \to M$.)

\begin{theorem}\label{theorem_spnat_21}
Let $M$ be a commutative monoid generated by the single element $a_0 \in M$.
Then there is a unique biadditive mapping $\btimes : M \times M \to M$ such that
$a_0 \btimes a_0 = a_0$.
Moreover, $\diamond$ is both associative and commutative and $a_0 \diamond a = a$ for all $a \in M$.
\end{theorem}

\proof
This is a special case of Theorem~\ref{theorem_bam_21}. However, 
the proof of the existence of the biadditive mapping $\btimes$ is rather simple and so
we give it here.

First note that if $\alpha_1,\,\alpha_2$ are endomorphisms of $M$ with $\alpha_1(a_0) = \alpha_2(a_0)$
then $\alpha_1 = \alpha_2$, since $\{ a \in M : \alpha_1(a) = \alpha_2(a) \}$ is a submonoid 
containing $a_0$, and is thus equal to $M$.
Note also that the mapping $\alpha_1 + \alpha_2 : M \to M$ given by 
$(\alpha_1 + \alpha_2)(a) = \alpha_1(a) + \alpha_2(a)$ for all $a \in M$ is an endomorphism of $M$
(since $M$ is commutative).

Let $M'$ be the set consisting of those elements
$a \in M$ for which there exists an endomorphism $\Lambda_a$ of $M$ such that
$\Lambda_a(a_0) = a$. Then $a_0 \in M'$, since we can take
$\Lambda_{a_0} = \id_M$. Moreover, $0 \in M'$ with $\Lambda_0 = 0$. 
Let $a_1,\,a_2 \in M'$, and hence  there exist endomorphisms $\Lambda_{a_1},\,\Lambda_{a_2}$
with $\Lambda_{a_1}(a_0) = a_1$ and $\Lambda_{a_2}(a_0) = a_2$. Here we can put
$\Lambda_{a_1 + a_2} = \Lambda_{a_1} + \Lambda_{a_2}$, since
$\Lambda_{a_1 + a_2}(a_0) = \Lambda_{a_1}(a_0) + \Lambda_{a_2}(a_0) = a_1 + a_2$, which
shows that $a_1 + a_2 \in M'$.
This implies that $M'$ is a submonoid of $M$ containing $a_0$
and therefore $M' = M$.
From the statement at the beginning of the proof the endomorphism $\Lambda_a$ 
with $\Lambda_a(a_0) = a$ is unique for each $a \in M$,
and from the above it then follows that $\Lambda_0 = 0$ and
$\Lambda_{a_1 + a_2} = \Lambda_{a_1} + \Lambda_{a_2}$ for all $a_1,\,a_2 \in M$.

Now define $\btimes : M \times M \to M$ by  $a_1 \btimes a_2 = \Lambda_{a_1}(a_2)$ for all 
$a_1,\,a_2 \in M$. Then the mapping $a_2 \mapsto a_1 \btimes a_2$ is a endomorphism for each $a_1 \in M$, 
since it is equal to $\Lambda_{a_1}$. Moreover, the mapping $a_1 \mapsto a_1 \btimes a_2 = \Lambda_{a_1}(a_2)$ 
is also an endomorphism for each $a_2 \in M$, since $\Lambda_0(a_2) = 0$ and
$\Lambda_{a_1 + b_1}(a_2) = \Lambda_{a_1}(a_2) + \Lambda_{b_1}(a_2)$
for all $a_1,\,b_1 \in M$.
Therefore the mapping $\btimes$ is biadditive and by definition
$a_0 \btimes a_0 = \Lambda_{a_0}(a_0) = a_0$. Moreover, it easily follows from
the statement at the beginning of the proof that it is the unique such bilinear mapping.

The  rest of the proof (namely that $\diamond$ is both associative and commutative and 
that $a_0 \diamond a = a$ for all $a \in M$)
is contained in the proof of Theorem~\ref{theorem_bam_21}.
\eop

Here is the application of Theorem~\ref{theorem_spnat_21} giving the existence of the multiplication
$\times$ only assuming that the principle of mathematical induction holds.

\begin{proposition}\label{prop_spnat_21}
Let $(X,f,x_0)$ be a minimal counting system with associated operation $+$.
Then there is a unique binary operation $\mtimes$ on $X$ satisfying
$\mathrm{(\mtimes_0)}$ and $\mathrm{(\mtimes_1)}$.
This operation is both associative and commutative and, 
in addition, the distributative law holds for $+$ and $\mtimes$.
\end{proposition}

\proof 
By Lemma~\ref{lemma_spnat_11} we can apply Theorem~\ref{theorem_spnat_21}
to the associated monoid $X$ and the element $f(x_0)$. This gives us a unique biadditive mapping 
$\times : X \times X \to X$ such that $f(x_0) \times f(x_0) = f(x_0)$. 
The operation $\times$ is associative and commutative with
$f(x_0) \times x = x$ for all $x \in X$. Moreover, the distributative law holds for $+$ and $\mtimes$,
since this is just a part of $\times$ being biadditive.
It is thus enough to show that $\times$ is the unique binary operation on $X$ satisfying
$\mathrm{(\mtimes_0)}$  and $\mathrm{(\mtimes_1)}$, and
clearly $\mathrm{(\mtimes_0)}$ holds.
Consider $x_1,\,x_2 \in X$; then by 
$\mathrm{(\star)}$ and since $f(x_0) \times x_2 = x_2$ 
\[ f(x_1) \mtimes x_2 = (f(x_0) + x_1) \mtimes x_2 = (f(x_0) \mtimes x_2) + (x_1 \mtimes x_2)
= x_2 + (x_1 \mtimes x_2) \]
and so $\mathrm{(\mtimes_1)}$ holds.
Hence $\mtimes$ satisfies $\mathrm{(\mtimes_0)}$  and $\mathrm{(\mtimes_1)}$.
If $\mtimes'$ is any other operation 
satisfying $\mathrm{(\mtimes_0)}$  and $\mathrm{(\mtimes_1)}$ then 
$X_0 = \{x \in X : \mbox{$x \mtimes y = x \mtimes' y$ for all $y \in X$} \}$
is $f$-invariant and contains $x_0$: 
$x_0 \mtimes y = x_0 = x_0 \mtimes' y$ by $\mathrm{(\mtimes_0)}$ and so $x_0 \in X_0$, and if
$x \in X_0$ then 
$f(x) \mtimes y = y + (x \mtimes y)  = y + (x \mtimes' y)  = f(x) \mtimes' y$
by $\mathrm{(\mtimes_1)}$, and hence $f(x) \in X_0$. Therefore $X_0 = X$, which shows that
${\mtimes'} = {\mtimes}$, i.e., $\mtimes$ is the unique operation satisfying
$\mathrm{(\mtimes_0)}$  and $\mathrm{(\mtimes_1)}$.
\eop

We now look again at the operation $+$, since it has
some further properties which are important when introducing the usual
order relation on the natural numbers. These properties are listed in Proposition~\ref{prop_spnat_31}.

\begin{proposition}\label{prop_spnat_31}
Let $(X,f,x_0)$ be a minimal counting system and let $+$ be the associated operation.

(1)\enskip
For all $x_1,\,x_2 \in X$ there exists an $x \in X$ such that either $x_1 = x + x_2$ or
$x_2 = x + x_1$.

(2)\enskip
If $x_0 \notin f(X)$ then $x_1 + x_2 \ne x_0$ for all $x_1,\,x_2 \in X$ with $x_2 \ne x_0$.

(3)\enskip
If $f$ is injective then  $+$ obeys the cancellation law.
\end{proposition}

\proof 
(1)\enskip
Let $X_0$ be the subset of $X$ consisting of those $x_1 \in X$ such that for each $x_2 \in X$ there exists 
$x \in X$ with either $x_1 = x + x_2$ or $x_2 = x + x_1$.
Thus $x_0 \in X_0$, since $x_2 = x_2 + x_0$ for all $x_2 \in X$.
Consider $x_1 \in X_0$ and let $x_2 \in X$; if $x_1 = x + x_2$ for some $x \in X$ then
$f(x_1) = f(x + x_2) = f(x) + x_2$. On the other hand, if 
$x_2 = x + x_1$ for some $x \in X$ then either $x = x_0$, in which case
\[ f(x_1) = f(x_0 + x_1) = f(x_2) = f(x_0 + x_2) = f(x_0) + x_2\;,\]
or $x \ne x_0$, and so by Lemma~\ref{lemma_spnat_21} below $x = f(x')$ for some $x' \in X$ and then
\[ x_2 = f(x') + x_1 = f(x' + x_1) = f(x_1 + x') = f(x_1) + x'
= x' + f(x_1)\;.\]
Therefore $f(x_1) \in X_0$. This shows $X_0$ is an $f$-invariant set containing $x_0$ and hence $X_0 = X$.

(2)\enskip
Let $x_2 \in X \setminus \{x_0\}$ and put
$X_0 = \{ x_1 \in X : x_1 + x_2 \ne x_0 \}$. 
Now the set $X_0$ is $f$-invariant, and in fact
$f(x_1) + x_2 = f(x_1 + x_2) \ne x_0$ for all $x_1 \in X$, since $x_0 \notin f(X)$.
Also $x_0 \in X_0$, since $x_0 + x_2 = x_2 \ne x_0$ and thus $X_0 = X$.
Hence $x_1 + x_2 \ne x_0$ for all $x_1 \in X$.

(3)\enskip This is part of Theorem~\ref{theorem_spnat_11}.
\eop

The above proof made use of the following fact, which also shows that
a minimal counting system is a Dedekind system if and only if the mapping 
is injective but not surjective.

\begin{lemma}\label{lemma_spnat_21}
If $(X,f,x_0)$ is minimal then
for each $x \in X \setminus \{x_0\}$ there exists an $x' \in X$ such that
$x = f(x')$. In particular, $f$ is surjective if and only if $x_0 \in  f(X)$.
\end{lemma}

\proof 
The set $\{x_0\} \cup f(X)$ is trivially $f$-invariant and contains $x_0$ and hence $X = \{x_0\} \cup f(X)$. 
\eop

If $(X,f,x_0)$ is a Dedekind system then it easily follows from Proposition~\ref{prop_spnat_31}
that for each $x_1,\,x_2 \in X$ exactly one 
of the following statements hold: (1) $x_1 = x_2$, 
(2) there exists a unique $x \in X \setminus \{x_0\}$ such that $x_1 = x + x_2$, and
(3) there exists a unique $x \in X \setminus \{x_0\}$ such that $x_2 = x + x_1$.
This is the trichotomy needed to define the usual order on $X$.

What if the principle of mathematical induction holds but one of the other Peano axioms does not, and so
either $x_0 \in f(X)$ or $f$ is not injective? These cases are dealt with in, for example, in
Section~7 of Preston \cite{preston},
but this is more-or-less what happens. In both cases $X$ is finite.
If $x_0 \in f(X)$ then $f$ is a bijection and the picture looks like:

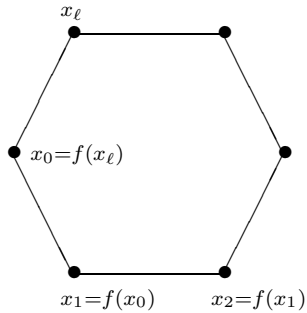
\setlength{\graphicthick}{0.1mm}
\setlength{\graphicmid}{0.1mm}
\setlength{\graphicthin}{0.1mm}

\begin{center}
\setlength{\unitlength}{0.8mm}
\begin{picture}(140,50)

\linethickness{\graphicthick}

\linethickness{\graphicthin}

\linethickness{\graphicmid}
\mythicklines

\put(49,24){$\bullet$}
\sput{53}{24}{x_0 = f(x_\ell)}

\sput{58}{0}{x_1 = f(x_0)}
\put(59,4){$\bullet$}

\sput{83}{0}{x_2 = f(x_1)}
\put(84,4){$\bullet$}

\put(59,44){$\bullet$}
\sput{58}{48}{x_\ell}

\put(84,44){$\bullet$}
\put(94,24){$\bullet$}

\put(50,25){\line(1,2){10}}
\put(50,25){\line(1,-2){10}}
\put(60,45){\line(1,0){25}}
\put(60,5){\line(1,0){25}}
\put(85,45){\line(1,-2){10}}
\put(85,5){\line(1,2){10}}

\end{picture}

\end{center}

Thus here $+$ and $\times$ are really nothing but addition and multiplication modulo $n$ with $n$ the 
cardinality of the finite set $X$
and, as stated in Theorem~\ref{theorem_spnat_11}, in this case the associated monoid
is an abelian group.
If $f$ is not injective then the picture is the following:

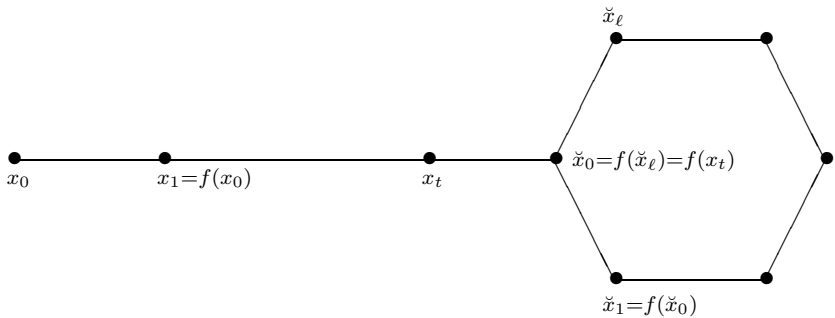
\begin{center}
\setlength{\unitlength}{0.8mm}
\begin{picture}(140,50)

\linethickness{\graphicthick}

\linethickness{\graphicthin}

\linethickness{\graphicmid}
\mythicklines

\put(5,25){\line(1,0){90}}

\put(95,25){\line(1,2){10}}
\put(95,25){\line(1,-2){10}}
\put(105,45){\line(1,0){25}}
\put(105,5){\line(1,0){25}}
\put(130,45){\line(1,-2){10}}
\put(130,5){\line(1,2){10}}

\sput{4}{21}{x_0}
\put(4,24){$\bullet$}

\sput{29}{21}{x_1 = f(x_0)}
\put(29,24){$\bullet$}

\sput{73}{21}{x_t}
\put(73,24){$\bullet$}

\sput{98}{24}{\breve{x}_0 = f(\breve{x}_\ell) = f(x_t)}
\put(94,24){$\bullet$}

\sput{103}{0}{\breve{x}_1 = f(\breve{x}_0)}
\put(104,4){$\bullet$}

\sput{103}{48}{\breve{x}_\ell}
\put(104,44){$\bullet$}

\put(129,4){$\bullet$}
\put(129,44){$\bullet$}
\put(139,24){$\bullet$}

\end{picture}

\end{center}

Even if it is surprising that an addition and a multiplication exist in this case it is a simple enough matter
to explicitly compute how they operate.

%%% Local Variables: 
%%% mode: latex
%%% TeX-master: "sums"
%%% End: 

\startsection{Minimal counting systems}

\label{mcs}

In the present section we do the following:

\begin{evlist}{18pt}{6pt}
\item[--]
Start by describing the more general set-up to be used in the remainder of these notes.
In this the single mapping occurring in a counting system will be replaced by a family of commuting mappings.

\item[--]
Formulate and prove the main result concerning the existence of an operation which is the analogue of
the addition on the natural numbers.

\end{evlist}

In Section~\ref{spnat} a counting system involved just a single mapping
$f : X \to X$. From now on we consider the more general situation in which
$f$ is replaced by a family of commuting mappings.
We therefore work with  triples $(X,\Fam{f},x_0)$ consisting of a set $X$, 
a non-empty family $\Fam{f} = \{f_s\}_{s \in S}$ of commuting mappings of $X$ into itself
and an element $x_0 \in X$. To be more precise about the family $\Fam{f}$,
we mean there is a non-empty index set $S$, for each $s \in S$ there is a mapping $f_s : X \to X$ and 
$f_s \circ f_t = f_t \circ f_s$ holds for all $s,\, t \in S$. 
Any such triple will be called an \definition{$S$-indexed counting system} or just a 
\definition{counting system} if $S$ can be determined from the context.
In the previous section $S$ consisted of a single element and $\Fam{f} = \{f\}$.

A subset $X'$ of $X$ is said to be \definition{$\Fam{f}$-invariant} if $f_s(X') \subset X'$ for all $s \in S$
and the counting system $(X,\Fam{f},x_0)$ is said to be \definition{minimal} if
the only $\Fam{f}$-invariant subset of $X$ containing 
$x_0$ is $X$ itself.

As a simple example of a counting system in which $S$ has more than one element consider the case which 
corresponds to the integers. Here there is a set $X$ (to be thought of as $\Int$) and a bijection 
$f_+ : X \to X$
(the successor operation). Put $f_- = f_+^{-1}$ (so $f_-$ is the predecessor operation)
and let $x_0 \in X$ (to be thought of as the integer $0$). Thus $S = \{+,-\}$, and 
$\Fam{f} = \{f_+,f_-\}$ is a family of commuting mappings, since $f_+ \circ f_- = f_- \circ f_+ = \id_X$. 
In the intended interpretation $(X,\Fam{f},x_0)$ will be a minimal counting system: The only subset of the 
integers containing $0$ which is invariant under both the successor and the predecessor operations is $\Int$
itself.

If $(X,\Fam{f},x_0)$ is an $S$-indexed counting system with $\Fam{f} = \{f_s\}_{s \in S}$ and
$(Y,\Fam{g},y_0)$ is a $T$-indexed counting system with $\Fam{g} = \{g_t\}_{t \in T}$ then
$(X \times Y,\Fam{f} \times \Fam{g},(x_0,y_0))$ is an $S \times T$-indexed counting system with the family
$\Fam{f} \times \Fam{g} = \{f_s \times g_t\}_{(s,t) \in S \times T}$, where
$(f_s \times g_t)(x,y) = (f_s(x),g_t(y))$ for all $x \in X$, $y \in Y$.
If $(X,\Fam{f},x_0)$ and $(Y,\Fam{g},y_0)$ are both minimal
then it is easy to see that $(X \times Y,\Fam{f} \times \Fam{g},(x_0,y_0))$ is also minimal.

The index set $S$ is now considered to be fixed, and so all counting systems are $S$-indexed counting
systems with this set $S$.
It is convenient to employ the following convention: If $(X,\Fam{f},x_0)$ is a counting system then it is 
assumed that the mapping in the family $\Fam{f}$ indexed by $s$ is always denoted by $f_s$, i.e., that
$\Fam{f} = \{f_s\}_{s \in S}$. Moreover, for each $s \in S$ the element $f_s(x_0)$ of $X$ will be denoted 
by $x_s$. (Thus if $(Y,\Fam{g},y_0)$ is a further counting system then $\Fam{g} = \{g_s\}_{s\in S}$  and $y_s$
denotes the element $g_s(y_0)$ of $Y$.)

Here is the main result concerning the existence of an operation $+$ which is the analogue of
the addition on the natural numbers:

\begin{theorem}\label{theorem_mcs_11}
Let $(X,\Fam{f},x_0)$ be a minimal counting system. Then there exists a unique binary operation $+$ such that
$(X,+,x_0)$ is a commutative monoid with
\begin{evlist}{12pt}{6pt}
\item[$\mathrm{(\star)}$]
$\ f_s(x) = x_s + x$ for all $x \in X$, $s \in S$.
\end{evlist}
Moreover, 
$X$ is an abelian group if and only if each mapping in $\Fam{f}$ is a bijection, and it 
obeys the cancellation law 
if and only if each mapping in $\Fam{f}$ is injective. 
\end{theorem}

\proof
This will follow from Theorem~\ref{theorem_mcs_21}. \eop

As in Section~\ref{spnat}
the operation $+$ in Theorem~\ref{theorem_mcs_11} will be called the
\definition{operation associated with $(X,\Fam{f},x_0)$}. 
It makes
$X$ a commutative monoid with zero $x_0$ and as such it is uniquely determined by the requirement that
$f_s(x) = x_s + x$ for all $x \in X$, $s \in S$. The monoid $X$ will also be referred to as
the \definition{associated monoid}.

The operation $+$ has the properties which correspond to 
$\mathrm{(+_0)}$ and $\mathrm{(+_1)}$ for the natural numbers:

\begin{proposition}\label{prop_mcs_11}
The operation $+$ associated with the minimal counting system $(X,\Fam{f},x_0)$ is  
the unique binary operation on $X$ such that
\begin{evlist}{12pt}{6pt}
\item[$\mathrm{(+_0)}$]
$\ x_0 + x = x$ for all $x \in X$,

\item[$\mathrm{(+_1)}$]
$\ f_s(x_1) + x_2 = f_s(x_1 + x_2)$ for all $x_1,\,x_2 \in X$ and each $s \in S$.
\end{evlist}
\end{proposition}

\proof 
This is the same as the proof of Proposition~\ref{prop_spnat_11}.
Clearly $\mathrm{(+_0)}$ holds, and
by $\mathrm{(\star)}$ it follows that
$f_s(x_1) + x_2 = (x_s + x_1) + x_2 = x_s + (x_1 + x_2) = f_s(x_1 + x_2)$
for all $x_1,\,x_2 \in X$ and each $s \in S$, i.e.,
$\mathrm{(+_1)}$ also holds.
Let $+'$ also satisfy $\mathrm{(+_0)}$ and 
$\mathrm{(+_1)}$ and put
$E = \{ y \in X : x +' y = x + y\ \mbox{for all $x \in X$} \}$.
Then $x_0 \in E$, since
$x_0 + x  = x = x_0 +' x$ for all $x \in X$, and if $y \in E$ and $s \in S$ then
$f_s(y) \in E$, since
$f_s(y) + x = f_s(y + x) = f_s(y +' x) = f_s(y) +' x$ for all $x \in X$.
Thus $E$ is an $\Fam{f}$-invariant subset of $X$ containing $x_0$, and hence $E = X$, 
which implies that ${+'} = {+}$. \eop

Consider the example which corresponds to the integers. Here we have  a
set $X$ and a bijection $f_+ : X \to X$. Put $f_- = f_+^{-1}$;
thus $\Fam{f} = \{f_+,f_-\}$ is a family of commuting mappings with $S = \{+,-\}$.
Let $x_0 \in X$ and assume the counting system $(X,\Fam{f},x_0)$ is minimal.
Then by Theorem~\ref{theorem_mcs_11} there exists a unique binary operation $+$ on $X$ such that
$(X,\Fam{f},x_0)$ is an abelian group and
$f_+(x) = x_+ + x$ and $f_-(x) = x_- + x$ for all $x \in X$,
where $x_+ = f_+(x_0)$ and $x_- = f_-(x_0)$.

\bigskip
As before denote the set of all mappings of a set $X$ into itself by $\Self{X}$;  
we thus have the monoid $(\Self{X},\circ,\id_X)$.
It turns out that there is a one-to-one correspondence between operations $+$
on $X$ for which $(X,+,x_0)$ a commutative monoid and certain commutative submonoids of
$\Self{X}$. This correspondence is essentially given by the analogue for monoids of Cayley's theorem for 
groups. The operation $+$ in Theorem~\ref{theorem_mcs_11} will be obtained 
via the corresponding submonoid of $\Self{X}$.

Let $M$ be a monoid. 
For each subset $N$ of $M$ there is a least submonoid of $M$ containing $N$ (namely the intersection
of all such submonoids) which will be denoted by $\lsm{N}$.
A subset $N$ of a $M$ is \definition{commutative} if 
$a_1 \bullet a_2 = a_2 \bullet a_1$
for all $a_1,\,a_2 \in N$. In particular, a submonoid of $M$ being \definition{commutative} 
means it is commutative as a subset of $M$. 

\begin{lemma}\label{lemma_mcs_11}
If $N \subset M$ is commutative then $\lsm{N}$ is a commutative submonoid.
\end{lemma}

\proof 
For each $a \in M$ the set $U_a = \{\,b \in M : \mbox{$b \bullet a = a \bullet b$}\, \}$ 
is a submonoid, since $e \bullet a = a = a \bullet e$ and if 
$a_1,\,a_2 \in U_a$ then 
\[ (a_1 \bullet a_2) \bullet a = a_1 \bullet a_2 \bullet a = a_1 \bullet a \bullet a_2 
= a \bullet a_1 \bullet a_2 = a \bullet (a_1 \bullet a_2)\;. \] 
Now $N \subset U_a$ for all $a \in N$, since $N$ is a commutative subset, and thus
$\lsm{N} \subset U_a$ for all $a \in N$, i.e.,
$b \bullet a = a \bullet b$ for all $a \in N$, $b \in \lsm{N}$.
But this also says that $N \subset U_b$ for each $b \in \lsm{N}$, which implies that
$\lsm{N} \subset U_b$ for each $b \in \lsm{N}$, and hence shows that
$b \bullet a = a \bullet b$ for all $a,\, b \in \lsm{N}$. \eop

For a counting system $(X,\Fam{f},x_0)$
the least submonoid of $\Self{X}$ containing $f_s$ for each  $s \in S$ will be denoted
by $M_\Fam{f}$.
By Lemma~\ref{lemma_mcs_11} $M_\Fam{f}$ is commutative, since
$\Fam{f}$ is a family of commuting mappings and hence
$\{ u \in \Self{X} : \mbox{$u = f_s$ for some $s \in S$} \}$ is a commutative subset of
$\Self{X}$. It is this commutative submonoid $M_\Fam{f}$ which will be used to obtain the operation $+$.

Let $\Phi_{x_0} : \Self{X} \to X$ be the evaluation mapping given by $\Phi_{x_0}(u) = u(x_0)$
for each $u \in \Self{X}$. The restriction of this mapping to any subset of $\Self{X}$ will also be
denoted by $\Phi_{x_0}$. Theorem~\ref{theorem_mcs_11} is really a special case of the next result:

\begin{theorem}\label{theorem_mcs_21}
If the counting system $(X,\Fam{f},x_0)$ is minimal
then there exists a unique binary operation $+$ on $X$ such that
$(X,+,x_0)$ is a commutative monoid and
$\Phi_{x_0} : M_\Fam{f} \to X$ is an isomorphism of monoids.
Moreover $u(x) = u(x_0) + x$ for all $x \in X$, $u \in M_\Fam{f}$. In particular
$\mathrm{(\star)}$ then holds, i.e., $f_s(x) = x_s + x$ for all $x \in X$, $s \in S$.
Furthermore, $X$ is an abelian group if and only if each mapping in $M_\Fam{f}$ is a 
bijection, and it obeys the cancellation law if and only if each mapping in $M_\Fam{f}$ is injective. 
\end{theorem}

\proof 
This is given below.
\eop

\textit{Proof of Theorem~\ref{theorem_mcs_11}\ }
Everything follows immediately from Theorem~\ref{theorem_mcs_21}, except for the
final statements and for these it is only necessary to note that
the sets $\{ u \in \Self{X} : \mbox{$u$ is injective} \}$ 
and $\{ u \in \Self{X} : \mbox{$u$ is a bijection} \}$ are both submonoids of $\Self{X}$. Hence
each mapping in $M_\Fam{f}$ is injective (resp.\ bijective)
if and only if each mapping in the family $\Fam{f}$ is injective (resp.\ bijective). \eop

Now to the proof of Theorem~\ref{theorem_mcs_21}.
In what follows consider the set $X$ and $x_0 \in X$ to be fixed; for the moment
the family $\Fam{f}$ is not involved.
An important property of the mapping $\Phi_{x_0} : \Self{X} \to X$ is that
\begin{evlist}{26pt}{6pt}
\item[$\mathrm{(\sharp)}$]
$\ u(\Phi_{x_0}(v)) = \Phi_{x_0}(u \circ v)\,$ for all $u,\,v \in \Self{X}$, 
\end{evlist}
since $u(v(x_0)) = (u \circ v)(x_0)$.

In the first few results below it is not necessary to assume the monoid operation is commutative, 
and as long as this is the case the operation will be denoted by $\bullet$ instead of $+$.
For each binary operation $\bullet$ on $X$ such that $(X,\bullet,x_0)$ is a monoid 
define a mapping $\Psi_\bullet : X \to \Self{X}$ by letting
$\Psi_\bullet(x)(x') = x \bullet x'$ for all $x,\,x' \in X$; put $M_\bullet = \Psi_\bullet(X)$.
The following lemma is the analogue of Cayley's theorem:

\begin{lemma}\label{lemma_mcs_21}
The mapping $\Psi_\bullet : X \to \Self{X}$ is an injective homomorphism, which implies that $M_\bullet$ 
is a submonoid of $\,\Self{X}$ and $\Psi_\bullet : X \to M_\bullet$ is an isomorphism.
Moreover, the inverse of the isomorphism $\Psi_\bullet$ is the
the mapping $\Phi_{x_0} : M_\bullet \to X$, and so in particular $\Phi_{x_0}$ is a bijection. 
\end{lemma}

\proof
The mapping $\Psi_\bullet$ is a homomorphism since if $x_1,\,x_2 \in X$ then 
\begin{eqnarray*}
\Psi_\bullet(x_1 \bullet x_2)(x) = (x_1 \bullet x_2) \bullet x 
&=& x_1 \bullet (x_2 \bullet x) = \Psi_\bullet(x_1)(x_2 \bullet x) \\
&=& \Psi_\bullet(x_1)(\Psi_\bullet(x_2)(x)) = (\Psi_\bullet(x_1) \circ \Psi_\bullet(x_2))(x) 
\end{eqnarray*}
for all $x \in X$, i.e., $\Psi_\bullet(x_1 \bullet x_2) = \Psi_\bullet(x_1) \circ \Psi_\bullet(x_2)$, and 
$\Psi_\bullet(x_0)(x) = x_0 \bullet x = x$ for all $x \in X$, i.e., $\Psi_\bullet(x_0) = \id_X$.
Moreover, $\Psi_\bullet$ is injective, since
\[ x_1 = x_1\bullet x_0 = \Psi_\bullet(x_1)(x_0) = \Psi_\bullet(x_2)(x_0) = x_2 \bullet x_0 = x_2 \]
whenever $\Psi_\bullet(x_1) = \Psi_\bullet(x_2)$.
Hence $M_\bullet = \Psi_\bullet(X)$ is a submonoid of $\Self{X}$ 
and the mapping $\Psi_\bullet : X \to M_\bullet$ is an isomorphism. 
Now the inverse of $\Psi_\bullet$ is the mapping $\Phi_{x_0} : M_\bullet \to X$, since
$(\Phi_{x_0} \circ \Psi_\bullet)(x) = \Phi_{x_0}(\Psi_\bullet(x)) = x \bullet x_0 = x = \id_X(x)$
for each $x \in X$ and  thus $\Phi_{x_0}$ is the inverse of $\Psi_\bullet$. \eop

The following is the converse of the construction in the previous result:

\begin{lemma}\label{lemma_mcs_31}
Let $M$ be a submonoid of $\,\Self{X}$ such that the mapping $\Phi_{x_0} : M \to X$ is a bijection.
Then there exists a unique binary 
operation $\bullet$ on $X$ such that $(X,\bullet,x_0)$ is a monoid and 
$\Phi_{x_0}$ is an isomorphism.
Moreover, $u(x) = \Phi_{x_0}(u) \bullet x$ for all $u \in M$, $x \in X$.
\end{lemma}

\proof 
Since $\Phi_{x_0} : M \to X$ is a bijection we can define a binary operation $\bullet$ on $X$ by letting
$x_1 \bullet x_2 = \Phi_{x_0}(\Phi_{x_0}^{-1}(x_1) \circ \Phi_{x_0}^{-1}(x_2))$
for all $x_1,\,x_2 \in X$. Then
\[ \Phi_{x_0}(u_1) \bullet \Phi_{x_0}(u_2) = \Phi_{x_0}(u_1 \circ u_2)\] 
for all $u_1,\,u_2 \in M$ and $\bullet$ is the unique binary operation  with this property.
The operation $\bullet$ is associative since $\circ$ has this property: If $x_1,\,x_2,\,x_3 \in X$ and
$u_1,\,u_2,\,u_3 \in M$ are such that $x_j = \Phi_{x_0}(u_j)$ for $j = 1,\,2,\,3$ then
\begin{eqnarray*}
(x_1 \bullet x_2) \bullet x_3 &=& (\Phi_{x_0}(u_1) \bullet \Phi_{x_0}(u_2)) \bullet \Phi_{x_0}(u_3) 
= \Phi_{x_0}(u_1 \circ u_2) \bullet \Phi_{x_0}(u_3)\\
 &=& \Phi_{x_0}( (u_1 \circ u_2) \circ u_3) = \Phi_{x_0}( u_1 \circ (u_2 \circ u_3))\\
&=& \Phi_{x_0}(u_1) \bullet \Phi_{x_0}(u_2 \circ u_3)
= \Phi_{x_0}(u_1) \bullet (\Phi_{x_0}(u_2) \bullet \Phi_{x_0}(u_3))\\ 
&=& x_1 \bullet (x_2 \bullet x_3) \;.
\end{eqnarray*}
Moreover, $x_0$ is the unit  for $\bullet$, since
if $x \in X$ and $u \in M$ is such that $x = \Phi_{x_0}(u)$ then
$ x \bullet x_0 = \Phi_{x_0}(u) \bullet \Phi_{x_0}(\id_X) = \Phi_{x_0}(u \circ \id_X)
= \Phi_{x_0}(u) =  x$, and in the same way $x_0 \bullet x = x$.
Thus $(X,\bullet,x_0)$ is a monoid and $\Phi_{x_0}$ is an isomorphism, and $\bullet$ is uniquely
determined by these two requirements.
Finally, let $u \in M$ and $x \in X$. Since $\Phi_{x_0}$ is surjective there exists $v \in M$
with $x = \Phi_{x_0}(v)$ and then by 
$\mathrm{(\sharp)}$
$u(x) = u(\Phi_{x_0}(v)) = \Phi_{x_0}(u \circ v)  = \Phi_{x_0}(u) \bullet \Phi_{x_0}(v)  
= \Phi_{x_0}(u) \bullet x$. \eop

If $M$ is a submonoid of $\Self{X}$ for which $\Phi_{x_0} : M \to X$ is a bijection
then the monoids $M$ and $X$ are isomorphic.
Hence they have the same algebraic properties. For example,
$M$ is commutative if and only if $X$ is,
$M$ is a group if and only if $X$ is and
$M$ obeys the left cancellation law if and only if $X$ does.
(A monoid $M$ obeys the left cancellation law if
$a_1 = a_2$ holds whenever $a \bullet a_1 = a \bullet a_2$ for some $a \in M$).

\begin{lemma}\label{lemma_mcs_41}
Let $M$ be a submonoid of $\,\Self{X}$ for which $\Phi_{x_0} : M \to X$ is a bijection.
Then:

(1)\enskip
If $u \in M$ is a bijection then $u^{-1} \in M$.

(2)\enskip
The monoid $M$ is a group  if and only if each mapping in $M$ is a bijection. 

(3)\enskip
The monoid $M$ obeys the left cancellation law  if and only if each mapping in $M$ is injective. 
\end{lemma}

\proof
(1)\enskip
Let $u \in M$ be a bijection. Then $u^{-1}(x_0) \in X$ and $\Phi_{x_0}$ is surjective and so 
there exists $v \in M$ with $\Phi_{x_0}(v) = u^{-1}(x_0)$; thus by $\mathrm{(\sharp)}$
\[  \Phi_{x_0}(u \circ v) = u(\Phi_{x_0}(v)) = u(u^{-1}(x_0)) = x_0 = \Phi_{x_0}(\id_X) \]
and therefore $u \circ v = \id_X$, since $\Phi_{x_0}$ is injective. 
This shows that $u^{-1} = v \in M$.

(2)\enskip
Clearly $M$ is a group if and only if  each mapping $u \in M$ is a bijection and 
$u^{-1} \in M$. But (1) implies that $u^{-1} \in M$ holds automatically whenever $u \in M$ is a bijection.

(3)\enskip
Suppose  $M$ obeys the left cancellation law. Let
$u \in M$ and $x_1,\,x_2 \in X$ with $u(x_1) = u(x_2)$. Then 
there exist
$u_1,\,u_2 \in M$ with $\Phi_{x_0}(u_1) = x_1$ and $\Phi_{x_0}(u_2) = x_2$ (since $\Phi_{x_0}$ is surjective), 
and hence by
$\mathrm{(\sharp)}$
\[   \Phi_{x_0}(u\circ u_1) = u(\Phi_{x_0}(u_1)) = u(x_1) = u(x_2) = u(\Phi_{x_0}(u_2)) 
= \Phi_{x_0}(u\circ u_2)\;.\] 
It follows that $u \circ u_1 = u \circ u_2$ (since $\Phi_{x_0}$ is injective) and so 
$u_1 = u_2$. In particular $x_1 = x_2$, which shows that $u$ is injective.
The converse is immediate, since 
if $u \in M$ is injective and $u \circ u_1 = u \circ u_2$ then $u_1 = u_2$.
\eop

Now if the above results are to be applied to the proof of Theorem~\ref{theorem_mcs_21} then 
the mapping $\Phi_{x_0} : M_\Fam{f} \to X$ will have to be a bijection.
As will be seen below, this is the case and it follows from the fact that $(X,\Fam{f},x_0)$ is a minimal 
counting system.

If $N \subset \Self{X}$ then a subset $X'$ of $X$ is said to be \definition{$N$-invariant} if 
$u(X') \subset X'$ for each $u \in N$. 
For each subset $A$ of $X$ there is a least $N$-invariant subset of $X$ containing $A$, namely
the intersection of all such subsets 
(noting that $X$ itself is always an $N$-invariant subset containing $A$).
In particular, there is a least $N$-invariant subset of $X$ containing 
the element $x_0$.

\begin{lemma}\label{lemma_mcs_51}
If $M$ is a commutative submonoid of $\,\Self{X}$ then $\Phi_{x_0} : M \to X$
is a bijection if and only if the only $M$-invariant subset of $X$ containing $x_0$ is $X$ itself.
\end{lemma}

\proof 
Denote by $I_M(x_0)$ the least $M$-invariant subset of $X$ containing $x_0$. The lemma thus states that
$\Phi_{x_0} : M \to X$ is a bijection if and only if $X = I_M(x_0)$.
Put $X_0 = \Phi_{x_0}(M)$; then $x_0 = \Phi_{x_0}(\id_X) \in X_0$, and 
if $x = \Phi_{x_0}(v) \in X_0$ then $u(x) = u(v(x_0)) = \Phi_{x_0}(u \circ v) \in X_0$ for all $u \in M$;
hence $X_0$ is an $M$-invariant subset of $X$ containing $x_0$, and so $I_M(x_0) \subset X_0$.
But each element of $X_0$ has the form $v(x_0)$ for some $v \in M$ and so lies in $I_M(x_0)$, since 
$I_M(x_0)$ is $M$-invariant and contains $x_0$. Therefore $X_0 = I_M(x_0)$
and in particular $\Phi_{x_0}$ is surjective (i.e., $X_0 = X$) if and only if $X = I_M(x_0)$. 
(Note that this statement holds for an arbitrary submonoid $M$ of $\Self{X}$.)

Suppose next that $X = I_M(x_0)$ and let $u_1,\,u_2 \in M$ with $\Phi_{x_0}(u_1) = \Phi_{x_0}(u_2)$, i.e.,
with $u_1(x_0) = u_2(x_0)$. Thus the set $X_0 = \{ x \in X : u_1(x) = u_2(x) \}$ contains $x_0$, and
it is $M$-invariant, since if $u_1(x) = u_2(x)$ then for all $v \in M$
\[ u_1(v(x)) = v(u_1(x)) = v(u_2(x)) = u_2(v(x))\]
(and here of course we require $M$ to be commutative).
Hence $X_0 = X$, i.e., $u_1 = u_2$, which implies that $\Phi_{x_0}$ is injective. 

The above thus shows that
$\Phi_{x_0}$ is surjective if and only if $X = I_M(x_0)$, and also that $\Phi_{x_0}$ is injective
whenever  $X = I_M(x_0)$. Therefore $\Phi_{x_0} : M \to X$ is a bijection if and only if
$X = I_M(x_0)$.
In fact, the proof shows in addition that $\Phi_{x_0}$ is bijective if and only if it is surjective.
\eop

\begin{lemma}\label{lemma_mcs_61}
Let $N$ be any subset of $\,\Self{X}$; then a subset of $X$ is $N$-invariant if and only if it is
$\lsm{N}$-invariant.
\end{lemma}

\proof 
Let $X' \subset X$ be $N$-invariant and put
$U = \{ u \in \Self{X} : u(X') \subset X' \}$; then
$\id_X \in U$, and if $u_1,\,u_2 \in U$ then
$(u_1 \circ u_2)(X') = u_1(u_2(X')) \subset u_1(X') \subset X'$, i.e., $u_1 \circ u_2 \in U$,
which shows that $U$ is a submonoid. Moreover $N \subset U$, since
$X'$ is $N$-invariant, and hence $\lsm{N} \subset U$. This shows that $X'$ is $\lsm{N}$-invariant.
The converse holds trivially.
\eop

\textit{Proof of Theorem~\ref{theorem_mcs_21}\ }
Put $F = \{ u \in \Self{X} : \mbox{$u = f_s$ for some $s \in S$} \}$; then a subset of $X$ is
$\Fam{f}$-invariant if and only if it is $F$-invariant, and thus
by Lemma~\ref{lemma_mcs_61} if and only if it is $M_\Fam{f}$-invariant, 
since  by definition $M_\Fam{f} = \lsm{F}$.
But the only $\Fam{f}$-invariant subset of $X$ containing $x_0$ is $X$ itself, since 
$(X,\Fam{f},x_0)$ is minimal, and hence $X$
is the only $M_\Fam{f}$-invariant subset of $X$ containing $x_0$.
Therefore by Lemma~\ref{lemma_mcs_51} the mapping $\Phi_{x_0} : M_\Fam{f} \to X$ is a bijection.
Now Lemma~\ref{lemma_mcs_31} can be applied to obtain a unique binary 
operation $+$ on $X$ such that $(X,+,x_0)$ is a monoid and 
$\Phi_{x_0} : (M_\Fam{f},\circ,\id_X) \to (X,+,x_0)$ is an isomorphism.
This monoid is commutative since $M_\Fam{f}$ is.
Moreover, $u(x) = \Phi_{x_0}(u) + x = u(x_0) + x$ for all $u \in M_\Fam{f}$, $x \in X$.
The final statements follow from Lemma~\ref{lemma_mcs_41}. \eop

%%% Local Variables: 
%%% mode: latex
%%% TeX-master: "sums"
%%% End: 

\startsection{Existence of biadditive mappings}

\label{bam}

In this section we do the following:

\begin{evlist}{18pt}{6pt}
\item[--]
Formulate and prove a result concerning the existence of biadditive mappings defined
on commutative monoids.
\item[--]
Apply this result to obtain operations which 
generalise the multiplication on the natural numbers.
\end{evlist}

A family $\{a_t\}_{t \in T}$ of elements from a monoid $M$ is said to \definition{generate} $M$ if
the only submonoid containing $a_t$ for each $t \in T$ is $M$ itself.
The results given below can be applied to the monoid associated with a minimal counting system
because of the following simple fact:

\begin{lemma}\label{lemma_bam_31}
If $(X,\Fam{f},x_0)$ is a minimal counting system then the family $\{x_s\}_{s \in S}$
generates the associated monoid
(with as always $x_s = f_s(x_0)$ for each $s \in S$).
\end{lemma}

\proof 
Let $X'$ be any submonoid of $X$ containing $x_s$ for each $s \in S$. If $x \in X'$ and
$s \in S$ then by $\mathrm{(\star)}$
$f_s(x) = x_s + x \in X'$, and thus $X'$ is an $\Fam{f}$-invariant
subset of $X$ which contains $x_0$ (since $X'$ is a submonoid). Thus $X' = X$.
\eop

In all of what follows let $M$ and $N$ be commutative monoids and
let $\{a_s\}_{s \in S}$ be a family of elements which generates $M$. 
A mapping $\btimes : M \times M \to N$ is called \definition{biadditive} if
$a' \mapsto a \btimes a'$ and $a' \mapsto a' \btimes a$ are both homomorphisms from $M$ to $N$
for each $a \in M$.
The question to be considered here is:
Given a mapping $\sigma : S \times S \to N$, does there exist a biadditive mapping
$\btimes : M \times M \to N$ such that
$a_s \btimes a_t = \sigma(s,t)$ for all $s,\,t \in S$?

Suppose that $\btimes : M \times M \to N$ is a biadditive mapping and for each
$s \in S$ define $\lambda_s,\,\lambda'_s : M \to N$ by
$\lambda_s(a) = a_s \btimes a$ and $\lambda'_s(a) = a \btimes a_s$ for all $a \in M$.
Then $\lambda_s$ and $\lambda'_s$ are homomorphisms and
$\lambda_s(a_t) = a_s \btimes a_t = \lambda'_t(a_s)$ for all $s,\,t \in S$.
A necessary condition for the existence of a biadditive mapping 
$\btimes : M \times M \to N$ is thus that for each $s \in S$ there exist
homomorphisms $\lambda_s$ and $\lambda'_s$ such that
$\lambda_s(a_t) = \lambda'_t(a_s)$ for all $s,\,t \in S$. 
Theorem~\ref{theorem_bam_11} below states that this requirement is in fact sufficient.

Let us denote the set of homomorphisms from $M$ to $N$ by $\Hom(M,N)$ and for 
$\alpha_1,\,\alpha_2 \in \Hom(M,N)$ define  a mapping $\alpha_1 + \alpha_2 : M \to N$ by letting
\[ (\alpha_1 + \alpha_2)(a) = \alpha_1(a) + \alpha_2(a)\]
for each $a \in M$. Then it is easily checked that $\alpha_1 + \alpha_2 \in \Hom(M,N)$ 
(since the monoid $N$ is commutative).

\begin{lemma}\label{lemma_bam_11}
If $\alpha_1,\,\alpha_2 \in \Hom(M,N)$ are homomorphisms with $\alpha_1(a_s) = \alpha_2(a_s)$ for all 
$s \in S$ then 
$\alpha_1 = \alpha_2$.
\end{lemma}

\proof
This follows since $\{ a \in M : \alpha_1(a) = \alpha_2(a) \}$ is a submonoid of $M$.
\eop

\begin{lemma}\label{lemma_bam_21}
If $\btimes,\, \btimes' : M \times M \to N$ 
are biadditive with $a_s \btimes  a_t = a_s \btimes' a_t$
for all $s,\, t \in S$ then ${\btimes'} = {\btimes}$. 
\end{lemma}

\proof 
By Lemma~\ref{lemma_bam_11} $a_s \btimes' a' = a_s \btimes a'$ for all $a' \in M$, and thus
$a_s \in M_0$ for all $s \in S$, where
$M_0 = \{ a \in M : \mbox{$a \btimes' a' = a \btimes a'$ for all $a' \in M$ for each $s \in S$} \}$.
But $M_0$ is a submonoid and hence $M_0 = M$, which shows that ${\btimes'} = {\btimes}$. \eop

\begin{theorem}\label{theorem_bam_11}
Suppose for each $s \in S$ there exist $\lambda_s,\,\lambda'_s \in \Hom(M,N)$
such that $\lambda_s(a_t) = \lambda_t'(a_s)$ for all $s,\,t \in S$.
Then there exists a unique biadditive mapping $\btimes : M \times M \to N$ such that
$a_s \btimes a_t = \lambda_s(a_t)$ for all $s,\,t \in S$. 
\end{theorem}

\proof
Let $M'$ be the set consisting of those elements
$a \in M$ for which there exists a homomorphism $\Lambda_a \in \Hom(M,N)$ such that
$\Lambda_a(a_s) = \lambda'_s(a)$ for all $s \in S$. Then $a_t \in M'$ for each $t \in S$,
since $\lambda_t(a_s) = \lambda'_s(a_t)$ for all $s \in S$ and so we can take
$\Lambda_{a_t} = \lambda_t$. Moreover, 
$0(a_s) = 0 = \lambda'_s(0)$ for all $s \in S$ and thus $0 \in M'$ with $\Lambda_0 = 0$. 
Let $a_1,\,a_2 \in M'$, so there exist $\Lambda_{a_1},\,\Lambda_{a_2} \in \Hom(M,N)$
with $\Lambda_{a_1}(a_s) = \lambda'_s(a_1)$ and $\Lambda_{a_2}(a_s) = \lambda'_s(a_2)$ 
for all $s \in S$. Put $\Lambda_{a_1 + a_2} = \Lambda_{a_1} + \Lambda_{a_2}$; then
\[\Lambda_{a_1 + a_2}(a_s) = \Lambda_{a_1}(a_s) + \Lambda_{a_2}(a_s)
= \lambda'_s(a_1) + \lambda'_s(a_2) = \lambda'_s(a_1 + a_2) \]
for all $s \in S$, and hence $a_1 + a_2 \in M'$.
This shows that $M'$ is a submonoid of $M$ with $a_s \in M'$ for each $s \in S$
and therefore $M' = M$.

Now Lemma~\ref{lemma_bam_11} implies that the homomorphism $\Lambda_a$ 
with $\Lambda_a(a_s) = \lambda'_s(a)$ for all $s \in S$ is unique for each $a \in M$,
and from the above it follows that $\Lambda_0 = 0$ and
$\Lambda_{a_1 + a_2} = \Lambda_{a_1} + \Lambda_{a_2}$ for all $a_1,\,a_2 \in M$.

Now define $\btimes : M \times M \to N$ by putting $a_1 \btimes a_2 = \Lambda_{a_1}(a_2)$ for all 
$a_1,\,a_2 \in M$. Then the mapping $a_2 \mapsto a_1 \btimes a_2$ is a homomorphism for each $a_1 \in M$, 
since it is equal to $\Lambda_{a_1}$. Moreover, the mapping $a_1 \mapsto a_1 \btimes a_2 = \Lambda_{a_1}(a_2)$ 
is also an homomorphism for each $a_2 \in M$, since $\Lambda_0(a_2) = 0$ and
$\Lambda_{a_1 + b_1}(a_2) = \Lambda_{a_1}(a_2) + \Lambda_{b_1}(a_2)$
holds for all $a_1,\,b_1 \in M$.
Therefore the mapping $\btimes$ is biadditive and by definition
$a_s \btimes a_t = \Lambda_{a_s}(a_t) = \lambda'_t(a_s) = \lambda_s(a_t)$ for all $s,\,t \in S$.

Finally, the uniqueness follows immediately from Lemma~\ref{lemma_bam_21}. \eop

As already noted, the existence of the endomorphisms 
$\lambda_s$ and $\lambda'_s$ for each $s\in S$ is a necessary condition for the existence of 
$\btimes$. Thus Theorem~\ref{theorem_bam_11} guarantees the existence of a biadditive operation 
provided it can be specified appropriately on the set $(M_S \times M) \cup (M \times M_S)$,
where $M_S = \{ a \in M : \mbox{$a = a_s$ for some $s \in S$} \}$.

We next consider the case with $N = M$.
If $M$ is a monoid then an \definition{endomorphism} of $M$ is a homomorphism $\alpha : M \to M$.

\begin{theorem}\label{theorem_bam_21}
Let $\odot : S \times S \to S$ be a binary operation on $S$ and 
suppose that for each $s \in S$ there exist endomorphims $\lambda_s$ and $\lambda'_s$ of $M$ 
with $\lambda_s(a_t) = a_{s\odot t}$ and $\lambda'_s(a_t) = a_{t\odot s}$ for all $t \in S$. 
Then there exists a unique biadditive mapping $\btimes : M \times M \to M$ such that
\[ a_s \btimes a_t = a_{s \odot t}\] 
for all $s,\,t \in S$. 
Moreover, if $\odot$ is associative resp.\ commutative then so is $\btimes$, and if
$\# \odot s = s$ for all $s \in S$ for some $\# \in S$ then $a_\# \btimes a = a$ for all $a \in M$.
\end{theorem}

\proof
By assumption $\lambda_s(a_t) = a_{s\odot t} = \lambda'_t(a_s)$
for all $s,\,t \in S$ and therefore by Theorem~\ref{theorem_bam_11}
there exists a unique biadditive mapping $\btimes : M \times M \to M$ 
such that $a_s \btimes a_t = \lambda_s(a_t) = a_{s\odot t}$ for all $s,\,t \in S$.
Suppose $\odot$ is associative; then
\[
a_r \btimes (a_s \btimes a_t) = a_r \btimes a_{s \odot t} = a_{r \odot (s \odot t)}
= a_{(r \odot s) \odot t} = a_{r \odot s} \btimes a_t = (a_r \btimes a_s) \btimes a_t
\]
for all $r,\,s,\,t \in S$. Fix $r \in S$ and define $\btimes_r$ and $\btimes'_r$ by
$a_1 \btimes_r a_2 = a_r \btimes (a_1 \btimes a_2)$ and
$a_1 \btimes'_r a_2 = (a_r \btimes a_1) \btimes a_2$ for all $a_1,\,a_2 \in M$.
Then $\btimes_r$ and $\btimes'_r$ are both biadditive and
$a_s \btimes_r a_t = a_s \btimes'_r a_t$ for all $s,\,t \in S$, and thus by Lemma~\ref{lemma_bam_21}
${\btimes_r} = {\btimes'_r}$, i.e., 
$a_r \btimes (a_1 \btimes a_2) = (a_r \btimes a_1) \btimes a_2$ for all $a_1,\,a_2 \in M$.
This means that $a_r \in M_0$ for all $r \in S$, where
$M_0 = \{ a \in M : 
\mbox{$a \btimes (a_1 \btimes a_2) = (a \btimes a_1) \btimes a_2$ for all $a_1,\,a_2 \in M$} \}$.
But $M_0$ is a submonoid and hence $M_0 = M$, which shows that $\btimes$ is associative. 
A similar (but easier) argument shows that $\btimes$ is commutative whenever $\odot$ is.
The final statement follows from the fact that if
$a_* \btimes a_s = a_s$
for some $a_* \in M$ and all $s \in S$ then $a_* \btimes a = a$ for all $a \in M$:
The mapping $a \mapsto a_* \btimes a$ is an
endomorphism which agrees with $\id_M$ on the elements of the family $\{a_s\}_{s \in S}$, and thus 
by Lemma~\ref{lemma_bam_11} is equal to $\id_M$. \eop

We again look at the example corresponding to the integers. Here we have  a
set $X$ and a bijection $f_+ : X \to X$. Put $f_- = f_+^{-1}$;
thus $\Fam{f} = \{f_+,f_-\}$ is a family of commuting mappings with $S = \{+,-\}$.
Let $x_0 \in X$ and assume the counting system $(X,\Fam{f},x_0)$ is minimal.
Then by Theorem~\ref{theorem_mcs_11} there exists a unique binary operation $+$ on $X$ such that
$(X,\Fam{f},x_0)$ is an abelian group and
$f_+(x) = x_+ + x$ and $f_-(x) = x_- + x$ for all $x \in X$,
where $x_+ = f_+(x_0)$ and $x_- = f_-(x_0)$.
Denote the negative of $x \in X$ in this group by $-x$.
Then $x_- = -x_+$ and $x_+ = -x_-$.
Note that the mappings $\id_X$ and $-\id_X$ (with $-\id_X(x) = -x$ for all $x \in X$) are both endomorphisms
of $X$.

Consider the binary operation $\odot : S \times S \to S$ with
${+} \odot {+} = {-} \odot {-} = {+}$ and ${+} \odot {-} = {-} \odot {+} = {-}$. 
Then $\odot$ is clearly both associative and commutative and
$\id_X(x_t) = x_{{+}\odot t}$, $-\id_X(x_t) = x_{{-}\odot t}$,
$\id_X(x_t) = x_{t\odot {+}}$ and $-\id_X(x_t) = x_{t\odot {-}}$ for all $t \in S$.
Thus by Theorem~\ref{theorem_bam_21} there 
exists a unique biadditive mapping $\mtimes : X \times X \to X$ with
$x_s \mtimes x_t = x_{s \odot t}$ for all $s,\,t \in S$, i.e., such that
\[x_+ \mtimes x_+ = x_- \mtimes x_- = x_+
\ \ \mbox{and}\ \ x_- \mtimes x_+ = x_+ \mtimes x_- = x_-\;.\]
Moreover,  $\mtimes$ is associative and commutative and 
$x_+ \mtimes x = x$ for all $x \in X$, since $+ \odot t = t$ for all $t \in S$.

%%% Local Variables: 
%%% mode: latex
%%% TeX-master: "sums"
%%% End: 

\startsection{Counting systems and commutative monoids}

\label{monoids}

\medskip

Recall that if $(X,\Fam{f},x_0)$ is a counting system then it is assumed 
that $\Fam{f} = \{f_s\}_{s \in S}$ and that for each $s \in S$ the element $f_s(x_0)$ of $X$ is denoted 
by $x_s$. 
Combining Theorem~\ref{theorem_mcs_11} and Lemma~\ref{lemma_bam_31} results in
the following statement: 

\begin{evlist}{26pt}{6pt}
\item[$\longrightarrow$]
If $(X,\Fam{f},x_0)$ is a minimal counting system then 
there exists a unique binary operation $+$ such that
$(X,+,x_0)$ is a commutative monoid with

\begin{evlist}{36pt}{6pt}
\item[$\mathrm{(\star)}$]
$\ f_s(x) = x_s + x$ for all $x \in X$, $s \in S$.
\end{evlist}

Moreover, the family $\{x_s\}_{s \in S}$ generates the monoid $X$, meaning that
the only submonoid containing $a_s$ for each $s \in S$ is $X$ itself. 
\end{evlist}

This procedure of starting with minimal counting system and ending up with a commutative monoid
and a generating family can be reversed. We will show below that the following holds:

\begin{evlist}{26pt}{6pt}
\item[$\longleftarrow$]
If $M$ is a commutative monoid and $\{a_s\}_{s \in S}$ is a family of elements which generates 
the monoid $M$ then, defining $\tau_s : M \to M$ by 

\begin{evlist}{36pt}{6pt}
\item[\phantom{$\mathrm{(\star)}$}]
$\ \tau_s(a) = a_s + a$ for all $a \in M$,
\end{evlist}

results in a minimal counting system $(M,\Fam{\tau},0)$.
\end{evlist}

By the uniqueness in Theorem~\ref{theorem_mcs_11} it follows that each of these procedures is the
inverse of the other. Moreover, it will be shown that they each
respect the structure preserving mappings. 
(The structure preserving mappings for counting systems are defined below; for monoids they
are of course the homomorphisms.)
This correspondence allows us to obtain results about minimal counting systems from results 
about commutative monoids, which tend to be easier to deal with.

First note two useful facts about monoids which will be needed several times.

\begin{lemma}\label{lemma_monoids_11}
Let $M$ be a monoid and $\{a_t\}_{t \in T}$ be a family generating $M$.

(1)\enskip
Is $B$ is any subset of $M$ containing $e$
such that $a_t \bullet b \in B$ for all $b \in B$ and all $t \in T$ then $B = M$.

(2)\enskip
Let $M'$ be a further monoid and let
$\alpha : M \to M'$ be a mapping with $\alpha(e) = e$ such that 
$\alpha(a_t \bullet b) = \alpha(a_t) \bullet \alpha(b)$ for all $b \in M$ and all $t \in T$. 
Then $\alpha : M \to M'$ is a homomorphism.
\end{lemma}

\proof
(1)\enskip
Let $N = \{ a \in M : \mbox{$a \bullet b \in B$ for all $b \in B$} \}$; then clearly
$e \in N$ and if $a_1,\,a_2 \in N$ then
$(a_1 \bullet a_2) \bullet b = a_1 \bullet (a_2 \bullet b) \in B$ for all $b \in B$, i.e.,
$a_1 \bullet a_2 \in N$. Thus $N$ is a submonoid of $M$ and by assumption $a_t \in N$ for each
$t \in T$; hence $N = M$.
But $N \subset B$, since $e \in B$, and therefore $B = M$.  

(2)\enskip
Let $N = \{ a \in M : \mbox{$\alpha(a \bullet b) = \alpha(a) \bullet' \alpha(b)$ for all $b \in B$} \}$; then
$e \in N$, since 
$\alpha(e \bullet b) = \alpha(b) = e' \bullet' \alpha(b) = \alpha(e) \bullet' \alpha(b)$ for all
$b \in B$, and if  $a_1,\,a_2 \in N$ then
\begin{eqnarray*} \alpha((a_1\bullet a_2) \bullet b) 
  &=& \alpha(a_1\bullet (a_2 \bullet b))   = \alpha(a_1)\bullet' \alpha(a_2 \bullet b)) 
  = \alpha(a_1)\bullet' (\alpha(a_2) \bullet' \alpha(b)) \\
  &=& (\alpha(a_1)\bullet' \alpha(a_2)) \bullet' \alpha(b) 
  = \alpha(a_1\bullet a_2) \bullet' \alpha(b) 
\end{eqnarray*}
for all $b \in M$, i.e., $a_1 \bullet a_2 \in N$.
Thus $N$ is a submonoid of $M$ and by assumption $a_t \in N$ for each
$t \in T$; hence $N = M$. Therefore $\alpha$ is a homomorphism. \eop

Now consider an arbitrary commutative monoid $M$ as well as an arbitrary family
$\{a_s\}_{s \in S}$ of elements from
$M$. For each $s \in S$ let $\tau_s : M \to M$ be the mapping given by
$\tau_s(a) = a_s + a$ for all $a \in M$ (i.e., $\tau_s$ is defined so that
$\mathrm{(\star)}$ holds).
Then we have the counting system $(M,\Fam{\tau},0)$ with $\Fam{\tau} = \{\tau_s\}_{s \in S}$.
Note that $\tau_s(0) =  a_s + 0 = a_s$ for each $s \in S$, which is compatible with the
previous convention for the meaning of the family $\{a_s\}_{s\in S}$.
We call $(M,\Fam{\tau},0)$ the 
\definition{counting system associated with $M$ and the family $\{a_s\}_{s \in S}$}.

\begin{proposition}\label{prop_monoids_11}
The counting system $(M,\Fam{\tau},0)$ is minimal if and only if the family $\{a_s\}_{s \in S}$ generates $M$.
\end{proposition}

\proof 
If $(M,\Fam{\tau},0)$ is minimal then the proof of Lemma~\ref{lemma_bam_31}
shows that the family $\{a_s\}_{s \in S}$ generates $M$ (since here 
$\mathrm{(\star)}$ holds by definition).
Thus suppose conversely that $\{a_s\}_{s \in S}$ generates $M$
and let $M'$ be an $\Fam{\tau}$-invariant subset of $M$ containing $0$. 
Then $0 \in M'$ and $a_s + a' = \tau_s(a') \in M'$ for all $a' \in M'$, $s \in S$.
Hence by Lemma~\ref{lemma_monoids_11}~(1) $M' = M$, and this shows that $(M,\Fam{\tau},0)$ is minimal. \eop

There is thus a one-to-one correspondence between minimal ($S$-indexed) counting systems
and pairs consisting of a commutative monoid and an ($S$-indexed) family of elements generating the monoid. 
More precisely, we have the following:

\begin{theorem}
(1)\enskip
Let $(X,+,x_0)$ be the monoid associated with the minimal counting system $(X,\Fam{f},x_0)$.
Then $(X,\Fam{f},x_0)$ is the counting system associated with the monoid $X$ and the family
$\{x_s\}_{s \in S}$.

(2)\enskip
Let $(M,\Fam{\tau},0)$ be the counting system associated with the commutative monoid $M$
and the generating family $\{a_s\}_{s \in S}$. Then $M$ is the monoid associated with the counting system
$(M,\Fam{\tau},0)$.

\end{theorem}

\proof 
(1)\enskip
By $\mathrm{(\star)}$ $f_s(x) = x_s + x$ for all $x \in X$, $s \in S$, and thus by
definition $(X,\Fam{f},x_0)$ is the counting system associated with the monoid $X$ and the family
$\{x_s\}_{s \in S}$.

(2)\enskip
The operation $+'$ associated with the minimal counting system
$(M,\Fam{\tau},0)$ is uniquely determined by the requirement that
$(M,+',0)$ is a commutative monoid with $\tau_s(a) = a_s +' a$ for all $a \in M$, $s \in S$.
But the monoid operation on $M$ also has these properties, and hence 
$M$ is the monoid associated with $(M,\Fam{\tau},0)$.
\eop

The above correspondence also carries over to mappings and 
to explain this we need the structure preserving mappings between counting systems.
If $(X,\Fam{f},x_0)$ and $(Y,\Fam{g},y_0)$ 
are counting systems then a mapping $\pi : X \to Y$ is said to be a
\definition{morphism} from $(X,\Fam{f},x_0)$ to $(Y,\Fam{g},y_0)$ if $\pi(x_0) = y_0$ and 
$g_s \circ \pi = \pi \circ f_s$ for each $s \in S$. 
It is clear that for each counting system $(X,\Fam{f},x_0)$ 
the identity mapping $\id_X$ is a morphism and it is easily checked that $\sigma \circ \pi$ is a morphism 
whenever $\pi : (X,\Fam{f},x_0) \to (Y,\Fam{g},y_0)$ and $\sigma : (Y,\Fam{g},y_0) \to (Z,\Fam{h},z_0)$ 
are morphisms. Moreover, if $\pi : (X,\Fam{f},x_0) \to (Y,\Fam{g},y_0)$ is a morphism then 
$\pi \circ \id_X = \pi = \id_Y \circ \pi$,
and if $\pi,\,\sigma$ and $\tau$ are morphisms for which the compositions are defined then
$(\tau \circ \sigma) \circ \pi = \tau \circ (\sigma \circ \pi)$.
This implies that
counting systems are the objects of a concrete category, whose morphisms are those just defined.

\begin{theorem}\label{theorem_monoids_11}
Let $(X,\Fam{f},x_0)$ and $(Y,\Fam{g},y_0)$ be minimal counting systems and let
$\pi : X \to Y$ be a mapping. Then $\pi : (X,\Fam{f},x_0) \to (Y,\Fam{g},y_0)$ is a morphism
if and only if $\pi : X \to Y$ is a homomorphism of the associated commutative monoids
with $\pi(x_s) = y_s$ for all $s \in S$. 
\end{theorem}

\proof 
Suppose first that $\pi : X \to Y$ is a homomorphism of the associated monoids with
$\pi(x_s) = y_s$ for all $s \in S$. If $x \in X$ and $s \in S$
then by $\mathrm{(\star)}$
\[ \pi(f_s(x)) = \pi(x_s + x) = \pi(x_s) + \pi(x) 
  = y_s + \pi(x) = g_s(\pi(x))\;, \]
which shows that $\pi \circ f_s = g_s \circ \pi$ for all $s \in S$. Thus
$\pi : (X,\Fam{f},x_0) \to (Y,\Fam{g},y_0)$ is a morphism (since also $\pi(x_0) = y_0$).

Now suppose that
$\pi : (X,\Fam{f},x_0) \to (Y,\Fam{g},y_0)$ is a morphism. Again $\pi(x_0) = y_0$ and
also $\pi(x_s) = \pi(f_s(x_0)) = g_s(\pi(x_0)) = g_s(y_0) = y_s$ for all $s \in S$.
Moreover
\[ \pi(x_s + x') = \pi(f_s(x')) = g_s(\pi(x')) = y_s + \pi(x')
                  =  \pi(x_s) +  \pi(x') \]
for all $x' \in X$ and all $s \in S$. Thus by Lemma~\ref{lemma_monoids_11}~(2) $\pi$ is a homomorphism, and
$\pi(x_s) = y_s$ for all $s \in S$. \eop

\textit{Important remark:}
Let us emphasise that in what follows all counting systems are
$S$-indexed counting systems for some fixed index set $S$.

A counting system $(X,\Fam{f},x_0)$ is said to be  \definition{initial} if for each counting system 
$(Y,\Fam{g},y_0)$ there exists a unique morphism from $(X,\Fam{f},x_0)$ to $(Y,\Fam{g},y_0)$.
This is the obvious generalisation of the definition given in Section~\ref{spnat} for the case of a single mapping,
where the recursion theorem was formulated as stating that a Dedekind system is initial.
In fact the converse also holds, i.e., any initial counting system is a Dedekind system
(Lawvere \cite{lawvere}).  
In Section~\ref{initial} we characterise initial counting systems
in the general case. This will be done by using Theorem~\ref{theorem_monoids_11}
to translate the problem into one involving commutative monoids, and in order to apply
Theorem~\ref{theorem_monoids_11} we first show in Lemma~\ref{lemma_monoids_61} below
that initial counting systems are minimal.

Let $(X,\Fam{f},x_0)$ be any counting system; then, since an arbitrary intersection of 
$\Fam{f}$-invariant
subsets of $X$ is again $\Fam{f}$-invariant and $X$ is itself an $\Fam{f}$-invariant subset containing $x_0$,
there is a least $\Fam{f}$-invariant subset of $X$ containing $x_0$ (namely the intersection
of all such subsets). Let us denote this subset by $X'$ and
for each $s \in S$
let $f'_s$ be the restriction of $f_s$ to the set $X'$ considered as a mapping in $\Self{X'}$. Put
$\Fam{f}' = \{f'_s\}_{s \in S}$.

\begin{lemma}\label{lemma_monoids_21}
The counting system $(X',\Fam{f}',x_0)$ is minimal.
\end{lemma}

\proof 
An $\Fam{f}'$-invariant subset $X'_0$ of $X'$ containing $x_0$ is also an $\Fam{f}$-invariant subset of $X$ 
containing $x_0$ and so
$X'\subset X'_0$. Thus $X'_0 = X'$, which implies that the only $\Fam{f}'$-invariant subset of $X'$ 
containing $x_0$ is $X'$ itself.
\eop

The counting system $(X',\Fam{f}',x_0)$ will be called the \definition{minimal core} of
$(X,\Fam{f},x_0)$.

\begin{lemma}\label{lemma_monoids_31}
If $(X,\Fam{f},x_0)$ is a minimal counting system then for each counting system
$(Y,\Fam{g},y_0)$ there at most one morphism $\pi : (X,\Fam{f},x_0) \to (Y,\Fam{g},y_0)$.
\end{lemma}

\proof 
If $\pi$ and $\pi'$ are both morphisms from $(X,\Fam{f},x_0)$ to $(Y,\Fam{g},y_0)$ then 
the set $\{ x \in X : \pi(x) = \pi'(x) \}$ is $\Fam{f}$-invariant and contains $x_0$
and it is thus equal to $X_0$, since $(X,\Fam{f},x_0)$ is  minimal.
Hence $\pi' = \pi$. \eop

A basic fact about initial objects (in any category) is that any two are isomorphic.
In the present situation
an \definition{isomorphism} is a morphism $\pi : (X,\Fam{f},x_0) \to (Y,\Fam{g},y_0)$ 
for which there exists a morphism $\sigma : (Y,\Fam{g},y_0) \to (X,\Fam{f},x_0)$ such that
$\sigma \circ \pi = \id_X$ and $\pi \circ \sigma = \id_Y$.
In this case $\sigma$ is uniquely determined by $\pi$ and is called the \definition{inverse} of $\pi$.
It is easily checked that
a morphism $\pi : (X,\Fam{f},x_0) \to (Y,\Fam{g},y_0)$ is an isomorphism if and only if
the mapping $\pi : X \to Y$ is a bijection, and that in this case the inverse morphism is the
inverse mapping $\pi^{-1} : Y \to X$.
Counting systems $(X,\Fam{f},x_0)$ and $(Y,\Fam{g},y_0)$ are said to be \definition{isomorphic} if there exists an
isomorphism $\pi : (X,\Fam{f},x_0) \to (Y,\Fam{g},y_0)$.

\begin{lemma}\label{lemma_monoids_41}
Any initial counting systems $(X,\Fam{f},x_0)$ and $(Y,\Fam{g},y_0)$ are isomorphic;
the unique morphism $\pi : (X,\Fam{f},x_0) \to (Y,\Fam{g},y_0)$ is an isomorphism. 
\end{lemma}

\proof 
Since $(Y,\Fam{g},y_0)$ is initial there exists a unique morphism $\sigma$ from $(Y,\Fam{g},y_0)$ to 
$(X,\Fam{f},x_0)$ and then $\sigma \circ \pi$ is a morphism from $(X,\Fam{f},x_0)$ to $(X,\Fam{f},x_0)$.
But $(X,\Fam{f},x_0)$ is initial and so there is a unique such morphism, which
is $\id_X$, and hence $\sigma \circ \pi = \id_X$. 
In the same way (reversing the roles of
$(X,\Fam{f},x_0)$ and $(Y,\Fam{g},y_0)$) it follows that $\pi \circ \sigma = \id_Y$ and therefore $\pi$ 
is an isomorphism.
\eop

\begin{lemma}\label{lemma_monoids_51}
An initial counting system is minimal.
\end{lemma}

\proof 
Let $(X,\Fam{f},x_0)$ be initial and let $(X',\Fam{f}',x_0)$ be its minimal core,
so $X'$ is the least $\Fam{f}$-invariant subset of $X$ containing $x_0$.
Now if $(Y,\Fam{g},y_0)$ is any counting system
then there exists a unique morphism $\pi : (X,\Fam{f},x_0) \to (Y,\Fam{g},y_0)$ and the
restriction $\pi_{|X'}$ of $\pi$ to $X'$ is a morphism
$\pi_{|X'} : (X',\Fam{f}',x_0) \to (Y,\Fam{g},y_0)$. Moreover, by Lemma~\ref{lemma_monoids_31}
it is the unique such morphism, which shows that 
$(X',\Fam{f}',x_0)$ is initial. But the inclusion mapping $i : X' \to X$
clearly defines a morphism from $(X',\Fam{f}',x_0)$ to $(X,\Fam{f},x_0)$, it is thus the unique morphism
and by Lemma~\ref{lemma_monoids_41} it is an isomorphism.
In particular $i$ is surjective, i.e., $X' = X$, which implies $(X,\Fam{f},x_0)$ is minimal.
\eop

We now come to the concept for commutative monoids which corresponds to a counting system being initial.
Let $M$ be a commutative monoid and $\{a_t\}_{t \in T}$ be a family of elements from $M$.
The monoid $M$ is said to be \definition{free with respect to $\{a_t\}_{t \in T}$} if
for each commutative monoid $N$ and each family $\{b_t\}_{t\in T}$ from $N$
there exists a unique homomorphism $\alpha : M \to N$ with $\alpha(a_t) = b_t$ for each $t \in T$.

\begin{lemma}\label{lemma_monoids_61}
If $M$ is free with respect to $\{a_t\}_{t \in T}$ then $\{a_t\}_{t \in T}$ generates $M$.
\end{lemma}

\proof 
This is the same as Lemma~\ref{lemma_monoids_51}:
Let $M$ be free with respect to $\{a_t\}_{t \in T}$, and let $M_0$ be the least
submonoid of $M$ containing $a_t$ for each $t \in T$. Then $M_0$ is also free
with respect to $\{a_t\}_{t \in T}$, and the result corresponding to Lemma~\ref{lemma_monoids_41} holds.
As in Lemma~\ref{lemma_monoids_51} it follows that the inclusion $i : M_0 \to M$ is surjective,
and hence $M_0 = M$. Thus $\{a_t\}_{t \in T}$ generates $M$.
\eop

By Lemma~\ref{lemma_monoids_61}
we can replace
`unique homomorphism' just by `homomorphism'
in the definition of $M$ being free with respect to $\{a_t\}_{t \in T}$, since
by Lemma~\ref{lemma_bam_11} 
any such homomorphism is automatically unique.

\begin{theorem}\label{theorem_monoids_21}
A minimal counting system $(X,\Fam{f},x_0)$ is initial if and only if the associated monoid
$X$ is free with respect to $\{x_s\}_{s \in S}$.
\end{theorem}

\proof 
Suppose first that $(X,\Fam{f},x_0)$ is initial and let
$N$ be a commutative monoid and $\{b_s\}_{s\in S}$ be a family of elements from $N$.
Let $N'$ be the least submonoid  of $N$ such that $b_s \in N'$ for all $s \in S$ and for each $s \in S$
let $\tau'_s \in \Self{N'}$ be given by $\tau'_s(b) = b_s \oplus b$ for all $b \in B'$; put 
$\Fam{\tau}' = \{\tau'_s\}_{s \in S}$.
By Proposition~\ref{prop_monoids_11} $(N',\Fam{\tau}',0)$ is then a minimal counting system, and
since $(X,\Fam{f},x_0)$ is initial there exists a morphism
$\pi : (X,\Fam{f},x_0) \to (N',\Fam{\tau}',0)$. Thus by Theorem~\ref{theorem_monoids_11}
$\pi$ is a homomorphism from the associated monoid $X$ to $N'$ with $\pi(x_s) = b_s$ for each $s \in S$ 
and so $\pi$,
regarded as a mapping from $X$ to $N$, is still a homomorphism.
Moreover, it is the unique such homomorphism with
$\pi(x_s) = b_s$ for each $s \in S$, since by 
Proposition~\ref{prop_monoids_11} $\{x_s\}_{s \in S}$ generates $X$.
This shows that $X$ is free with respect to $\{x_s\}_{s \in S}$.

Suppose conversely that the associated monoid $X$ is free with respect to $\{x_s\}_{s \in S}$.
Let $(Y,\Fam{g},y_0)$ be a counting system
and let $(Y',\Fam{g}',y_0)$ be its minimal core.
Then there exists a unique homomorphism $\alpha$ from $X$ to the associated monoid $Y'$ 
with $\alpha(x_s) = y_s$ for each $s \in S$. 
Therefore by Theorem~\ref{theorem_monoids_11}
$\alpha : (X,\Fam{f},x_0) \to (Y',\Fam{g}',y_0)$ is a morphism and 
so $\alpha$, regarded as a mapping from $X$ to $Y$, is a morphism
from $(X,\Fam{f},x_0)$ to $(Y',\Fam{g}',y_0)$.
By Lemma~\ref{lemma_monoids_31} it is the unique such morphism and hence $(X,\Fam{f},x_0)$ is initial.
\eop

We end the section by looking at the special case of a counting system with a single mapping
(i.e., the set-up we started with in Section~\ref{spnat}).
Let $(X,f,x_0)$ be such a counting system, so $f \in \Self{X}$. If
$(X,f,x_0)$ is minimal then by Proposition~\ref{prop_monoids_11}
the associated monoid is generated by the single
element $f(x_0)$, i.e., the only submonoid containing $f(x_0)$ is $X$ itself.
Conversely, if $M$ is a commutative monoid generated by a single element $a_0$
then, again by Proposition~\ref{prop_monoids_11},
the associated counting system $(M,\tau,0)$ is minimal, and here $\tau \in \Self{M}$
is the mapping given by $\tau(a) = a_0 + a$ for all $a \in M$.

A commutative monoid $M$ is said to be \definition{free with respect to $a_0 \in M$} if
for each commutative monoid $N$ and each $b \in N$
there exists a unique homomorphism $\alpha : M \to N$ with $\alpha(a) = b$.
This is just the case when the family consists of a single element
(with the braces being omitted).
If $M$ is free with respect to $a_0$ then by Lemma~\ref{lemma_monoids_61}
$M$ is generated by $a_0$.

\begin{proposition}\label{prop_monoids_21}
A commutative monoid $M$ generated by a single element $a_0$ is free with
respect to $a_0$ if and only if the associated counting system $(M,\tau,0)$ is initial. 
\end{proposition}

\proof 
This is a special case of Theorem~\ref{theorem_monoids_21}.
\eop

The recursion theorem and its converse state that a counting system with a single mapping is
initial if and only if it is a Dedekind system. (This will be proved in Section~\ref{initial}.)
Thus a commutative monoid $M$ generated by a single element $a_0$ is free with respect to $a_0$
if and only if the associated counting system $(M,\tau,0)$ is a Dedekind system.
By Proposition~\ref{prop_spnat_31} it easily follows that this is the case
if and only if the cancellation law holds in $M$ and $a_1 + a_2 = 0$ is only possible with
$a_1 = a_2 = 0$.

%%% Local Variables: 
%%% mode: latex
%%% TeX-master: "sums"
%%% End: 

\startsection{Free commutative monoids}

\label{free}

By Theorem~\ref{theorem_monoids_21}
a minimal counting system $(X,\Fam{f},x_0)$ is initial if and only if the associated monoid
$X$ is free with respect to the family $\{x_s\}_{s \in S}$.
As a preparation for characterising initial counting systems in the next section
we thus look here at the corresponding characterisation
of free commutative monoids.

Recall from the end of the previous section that
a commutative monoid $M$ is said to be \definition{free with respect to an element $a \in M$} if
for each commutative monoid $N$ and each $b \in N$
there exists a unique homomorphism $\alpha : M \to N$ with $\alpha(a) = b$.
By Lemma~\ref{lemma_monoids_61} a necessary condition for this to hold is that
$M$ be generated by the single element $a$.

A further concept that plays a role here is that of an internal direct sum.
Let $M$ be a commutative monoid and $\{M_t\}_{t \in T}$ be a family of submonoids of $M$. Then
$M$ is called the \definition{internal direct sum} of the family $\{M_t\}_{t \in T}$ if for
each commutative monoid $N$ and for each family $\{\alpha_t\}_{t \in T}$ with
$\alpha_t \in \Hom(M_t,N)$ for each $t \in T$ there exists a unique homomorphism
$\alpha \in \Hom(M,N)$ such that $\alpha_t$ is the restriction of $\alpha$ of $M_t$ for each $t \in T$.  

In what follows let $M$ be a commutative monoid and let $\{a_s\}_{s \in S}$ be a family of elements which generates 
$M$; let $M_s$ be the least monoid containing the element $a_s$ for each $s \in S$.

The main results in this section are Proposition~\ref{prop_free_21} and
Theorem~\ref{theorem_free_21}; the crucial technical result is 
Theorem~\ref{theorem_free_11}.
Note the following special case of Theorem~\ref{theorem_bam_11}:

\begin{proposition}\label{prop_free_11}
Suppose that for each $s \in S$ there exists $\delta_s \in \Hom(M,N)$
such that $\delta_s(a_t) = 0$ for all $t \ne s$.
Then there exists a unique biadditive mapping $\bigtriangleup : M \times M \to N$ such that
$a_s \bigtriangleup a_t = 0$ for all $t \ne s$ and $a_s \bigtriangleup a_s = \delta_s(a_s)$ for all $s \in S$.
\end{proposition}

\proof 
This follows from Theorem~\ref{theorem_bam_11} with $\lambda'_s = \lambda_s = \delta_s$ for all $s \in S$. 
\eop

\begin{theorem}\label{theorem_free_11}
Suppose there exists a biadditive mapping $\bigtriangleup : M \times M \to N$ with $a_s \bigtriangleup a_t = 0$
for all $s \ne t$. Then there exists a unique  $\alpha \in \Hom(M,N)$
such that $\alpha(a_s) = a_s \bigtriangleup a_s$ for all $s \in S$.
\end{theorem}

\proof Below. \eop

\begin{proposition}\label{prop_free_21}
The monoid $M$ is the internal direct sum
of $\{M_s\}_{s \in S}$ if and only if there exists
a unique biadditive mapping $\bigtriangleup : M \times M \to M$ such that
$a_s \bigtriangleup a_s = a_s$ for all $s \in S$ and $a_s \bigtriangleup a_t = 0$ for all $s \ne t$.
\end{proposition}

\proof
Suppose first that $M$ is the internal direct sum of the family $\{M_s\}_{s \in S}$.
Let $s \in S$; for $t \ne s$ let $\alpha'_t = 0$ as element of $\Hom(M_t,M)$
and let $\alpha'_s \in \Hom(M_s,M)$ be the inclusion mapping.
Then there exists 
$\delta_s \in \Hom(M,N)$ such that $\alpha'_t$ is the restriction of $\delta_t$ of $M_t$ for each $t \in S$.  
In particular, $\delta_s(a_s) = a_s$ and $\delta_s(a_t) = 0$ for all $t \ne s$.
Thus by Proposition~\ref{prop_free_11}
there exists a unique biadditive mapping $\bigtriangleup : M \times M \to N$ such that
$a_s \bigtriangleup a_s = \delta_s(a_s) = a_s$ for all $s \in S$ and
$a_s \bigtriangleup a_t = 0$ for all $t \ne s$.

Suppose conversely there exists a unique biadditive mapping $\bigtriangleup : M \times M \to M$ such that
$a_s \bigtriangleup a_s = a_s$ for all $s \in S$ and $a_s \bigtriangleup a_t = 0$ for all $s \ne t$;
let $N$ be a commutative monoid and $\{\alpha_s\}_{s \in S}$ a family with
$\alpha_s \in \Hom(M_s,N)$ for each $s \in S$.
Also for each $s \in S$ let $\delta_s$ be the endomorphism of $M$ with
$\delta_s(a) = a \bigtriangleup a_s$ for all $a \in M$. Then $\delta_s(a_s) = a_s$ and
$\delta_s(a_t) = 0$ for all $t \ne s$, from which it follows that $\delta_s(M) \subset M_s$ (since
$\{ a \in M : \delta_s(a) \in M_s \}$ is a submonoid  containing $a_t$ for each $t \in S$).
Define a homomorphism $\delta'_s : M \to N$ by
$\delta'_s(a) = \alpha_s(\delta_s(a))$ for all $a \in M$, thus $\delta'_s(a_s) = \alpha_s(a_s)$
and $\delta'_s(a_t) = \alpha_s(0) = 0$ for all $t \ne s$.
Therefore by Theorem~\ref{theorem_bam_11} there exists a unique
biadditive mapping $\bigtriangleup' : M \times M \to N$ such that
$a_s \bigtriangleup' a_s = \alpha_s(a_s)$ for all $s \in S$ and $a_s \bigtriangleup' a_t = 0$ for all $s \ne t$,
and hence by Theorem~\ref{theorem_free_11} there exists a unique homomorphism $\alpha \in \Hom(M,N)$
such that $\alpha(a_s) = \alpha_s(a_s)$ for all $s \in S$. It then follows that
$\alpha(a) = \alpha_s(a)$ for all $a \in M_s$ for each $s \in S$ and that it is the unique such homomorphism.
This shows that $M$ is the internal direct sum of the family $\{M_s\}_{s \in S}$.
\eop

\begin{proposition}\label{prop_free_31}
Suppose that for each $s \in S$ there exists an endomorphism $\delta_s$ of $M$
such that $\delta_s(a_s) = a_s$ and $\delta_s(a_t) = 0$ for all $t \ne s$.
Then the monoid $M$ is the internal direct sum
of the family $\{M_s\}_{s \in S}$.
\end{proposition}

\proof 
This follows directly from Propositions \ref{prop_free_11} and \ref{prop_free_21}
and Theorem~\ref{theorem_free_11}.
\eop

\begin{theorem}\label{theorem_free_21}
The following are equivalent:

(1)\enskip
The monoid $M$ is free with respect to $\{a_s\}_{s \in S}$.

(2)\enskip
The monoid $M_s$ is free with respect to $a_s$ for each $s \in S$ and 
$M$ is the internal direct sum of the family $\{M_s\}_{s \in S}$.

(3)\enskip For each commutative monoid $N$ and each family $\{b_s\}_{s\in S}$ from $N$
there exists a unique biadditive mapping $\bigtriangleup : M \times M \to N$ such that
$a_s \bigtriangleup a_s = b_s$ for all $s \in S$ and $a_s \bigtriangleup a_t = 0$ for all $s \ne t$.
\end{theorem}

\proof
(2) $\Rightarrow$ (1):
Let $N$ be a commutative monoid $N$ and let $\{b_t\}_{t\in T}$  a family from $N$.
Then, since $M_s$ is free with respect to $a_s$, there exists a homomorphism $\alpha_s \in \Hom(M_s,N)$
with $\alpha_s(a_s) = b_s$ for each $s \in S$.
Hence, since $M$ is the internal direct sum of the family $\{M_s\}_{s \in S}$, there exists 
$\alpha \in \Hom(M,N)$ with $\alpha_s$ the restriction of $\alpha$ to $M_s$ for each $s \in S$. Thus
$\alpha(a_s) = b_s$ for each $s \in S$ and by Lemma~\ref{lemma_bam_11} $\alpha$ is the unique such
homomorphism. This shows that $M$ is free with respect to $\{a_s\}_{s \in S}$.

(1) $\Rightarrow$ (3): 
There exists for each $s \in S$ a homomorphism $\delta_s \in \Hom(M,N)$ 
with $\delta_s(a_s) = b_s$ and $\delta_s(a_t) = 0$ for all $t \ne s$ and thus by 
Proposition~\ref{prop_free_11} there exists a unique
biadditive mapping $\bigtriangleup : M \times M \to N$ such that
$a_s \bigtriangleup a_s = b_s$ for all $s \in S$ and $a_s \bigtriangleup a_t = 0$ for all $s \ne t$.

(3) $\Rightarrow$ (2): 
It follows immediately from Proposition~\ref{prop_free_21} that
$M$ is the internal direct sum of the family $\{M_s\}_{s \in S}$, and so it remains to show that $M_s$
is free with respect to $a_s$ for each $s \in S$.
Fix $s \in S$ and let $N$ be a commutative monoid and $b \in N$;
put $b_s = b$ and for $t \ne s$ put $b_t = 0$.
There thus exists a biadditive mapping $\bigtriangleup : M \times M \to N$ such that
$a_r \bigtriangleup a_r = b_r$ for all $r \in S$ and $a_r \bigtriangleup a_t = 0$ for all $r \ne t$.
Define $\alpha \in \Hom(M_s,N)$ by $\alpha(a) = a \bigtriangleup a_s$ for each $a \in M_s$.
Then $\alpha(a_s) = b$ and by Lemma~\ref{lemma_bam_11} $\alpha$ is the unique such
homomorphism. This shows that $M_s$ is free with respect to $a_s$. \eop

By Proposition~\ref{prop_monoids_21} the monoid $M_s$ is free with respect to
$a_s$ if and only if the associated counting system $(M_s,\tau_s,0)$ is initial,
where $\tau_s : M_s \to M_s$ is given by $\tau_s(a) = a_s + a$ for all $a \in M_s$.
In Section~\ref{initial} we will see that this is the case if and only if $(M_s,\tau_s,0)$
is a Dedekind system.

\bigskip

\textit{Proof of Theorem~\ref{theorem_free_11}}
\ The homomorphism $\alpha \in \Hom(M,N)$ will be obtained by patching together homomorphisms defined on  
suitable submonoids of $M$. For each $E \subset S$ let $M_E$ be the least submonoid 
of $M$ containing the element 
$a_s$ for each $s \in E$. Thus $M_\varnothing = \{0\}$ and
if $F \subset E$ then clearly $M_F \subset M_E$.
Define a subset $\mathcal{G}$ of $\mathcal{P}(S)$ by 
\[\mathcal{G} = \{ E \subset S :
  \mbox{there exists $a_E \in M_E$ such that 
        $a_s \bigtriangleup a_s = a_s \bigtriangleup a_E$ for all $s \in E$} \}\;. \]

\begin{lemma}\label{lemma_free_11}
If $E \in \mathcal{G}$ then there is a unique homomorphism $\alpha_E \in \Hom(M_E,N)$
such that $\alpha_E(a_s) = a_s \bigtriangleup a_s$ for all $s \in E$.
Moreover, if $E,\,F \in \mathcal{G}$ with $F \subset E$ then $\alpha_F$ is the restriction of
$\alpha_E$ to $M_F$.
\end{lemma}

\proof
Since $E \in \mathcal{G}$ there there exists $a_E \in M_E$ such that 
$a_s \bigtriangleup a_s = a_s \bigtriangleup a_E$ for all $s \in E$. Thus
a homomorphism $\alpha_E \in \Hom(M_E,N)$ can be defined
by letting $\alpha_E(a) = a \bigtriangleup a_E$ for all $a \in M_E$ and then
$\alpha_E(a_s) = a_s \bigtriangleup a_E = a_s \bigtriangleup a_s$ for all $s \in E$.
Moreover, by Lemma~\ref{lemma_bam_11} $\alpha_E$ is the unique element of
$\Hom(M_E,N)$ with $\alpha_E(a_s) = a_s \bigtriangleup a_s$ for all $s \in E$ (even though
$a_E$ is not necessarily unique.)
The final statement also follows from Lemma~\ref{lemma_bam_11}. \eop

The following concept now plays an important role:
A subset $\mathcal{S}$ of $\mathcal{P}(S)$ is called an \definition{inductive system} if 
$\varnothing \in \mathcal{S}$ and $E \cup \{s\} \in \mathcal{S}$ for all $E \in \mathcal{S}$ and all $s \in S$.

\begin{lemma}\label{lemma_free_21}
The set $\mathcal{G}$ is an inductive system.
\end{lemma}

\proof
Clearly $\varnothing \in \mathcal{G}$ (with $a_\varnothing = 0$), and so consider $E \in \mathcal{G}$ and 
$t \in S$. We need to show that $E \cup \{t\} \in \mathcal{G}$ and thus it can be assumed that 
$t \in S \setminus E$. 
Then $a_t \bigtriangleup a = 0$ for all $a \in M_E$, since
$\{ a \in M : a_t \bigtriangleup a = 0 \}$ is a submonoid of $M$ and
$a_t \bigtriangleup a_s = 0$ for all $s \in E$.
Let $a_E \in M_E$ be such that $a_s \bigtriangleup a_s = a_s \bigtriangleup a_E$ for all $s \in E$ and put
$a_{E \cup \{t\}} = a_E + a_t$.
Therefore $a_{E \cup \{t\}} \in M_{E \cup \{t\}}$ and 
\[ a_s \bigtriangleup a_{E \cup \{t\}} = a_s \bigtriangleup (a_E + a_t)
   = (a_s \bigtriangleup a_E) + (a_s \bigtriangleup a_t) = a_s + a_s \]
for all $s \in E \cup \{t\}$, since
$a_s \bigtriangleup a_t = 0$ for all $s \in E$ and $a_t \bigtriangleup a_E = 0$.
Hence $E \cup \{t\} \in \mathcal{G}$. This shows that $\mathcal{G}$ is an inductive system.
\eop

\begin{lemma}\label{lemma_free_31}
If $\mathcal{S}$ is any inductive system then $M = \bigcup_{E \in \mathcal{S}} M_E$.
\end{lemma}

\proof
Put $M_0 = \bigcup_{E \in \mathcal{S}} M_E$; then $0 \in M_\varnothing \subset M_0$
and if $a \in M_0$ and $s \in S$ then $a \in M_E$ for some $E \in \mathcal{S}$ and hence
$a_s + a \in M_{E \cup \{s\}} \subset M_0$. 
Thus by Lemma~\ref{lemma_monoids_11}~(1) $M_0 = M$.
\eop

One last concept is needed: A subset $\mathcal{S}$ of $\mathcal{P}(S)$ is said to be \definition{directed} if
for all $E,\,F \in \mathcal{S}$ there exists $G \in \mathcal{S}$ with
$E \cup F \subset G$. Suppose now there exists a subset $\mathcal{G}'$ of $\mathcal{G}$ which
is both inductive and directed. Then the homomorphism $\alpha \in \Hom(M,N)$ can be defined as follows:
Let $a \in M$; by Lemma~\ref{lemma_free_31} there exists
$E \in \mathcal{G}'$ with $a \in E$. Moreover, if $a$ also lies in  some other
$F \in \mathcal{G}'$ then, since $\mathcal{G}'$ is directed, there exists $G \in \mathcal{G}'$ with 
$E \cup F \subset G$, and so by the final statement in 
Lemma~\ref{lemma_free_11} $\alpha_E(a) = \alpha_G(a) = \alpha_F(a)$.
This implies there is a unique mapping $\alpha : M \to N$ such that $\alpha(a) = \alpha_E(a)$ for all
$a \in M_E$, $E \in \mathcal{G}'$.

\begin{lemma}\label{lemma_free_41}
$\alpha$ is the unique homomorphism with
$\alpha(a_s) = a_s \bigtriangleup a_s$ for all $s \in S$. 
\end{lemma}

\proof 
Clearly $\alpha(0) = 0$. Let $a_1,\,a_2 \in M$; then, since $\mathcal{G}'$ is directed, there exists
$E \in \mathcal{G}'$ such that both $a_1$ and $a_2$ lie in $M_E$. Thus $a_1 + a_2 \in M_E$
and therefore
$\alpha(a_1 + a_2) = \alpha_E(a_1 + a_2) = \alpha_E(a_1) + \alpha_E(a_2) 
= \alpha(a_1) + \alpha(a_2)$.
This shows $\alpha$ is a homomorphism. Moreover,
if $s \in S$ and $E$ is any element of $\mathcal{G}'$ containing $a_s$ then
$\alpha(a_s) = \alpha_E(a_s) = a_s \bigtriangleup a_s$.
Finally, the uniqueness follows immediately from Lemma~\ref{lemma_bam_11}. \eop

The proof of Theorem~\ref{theorem_free_11} can thus be completed by exhibiting a
subset of $\mathcal{G}$ which is both inductive and directed. 
Note that an arbitrary intersection of inductive systems is again an inductive system,
and so there exists a least inductive system $\mathcal{F}$. In particular
$\mathcal{F} \subset \mathcal{G}$.

\begin{lemma}\label{lemma_free_51}
The least inductive system $\mathcal{F}$ is directed. In fact $A \cup B \in \mathcal{F}$
for all $A,\,B \in \mathcal{F}$.
\end{lemma}

\proof 
Consider $B \in \mathcal{F}$ to be fixed and let $\mathcal{S} = \{ A \subset S : A \cup B \in \mathcal{F} \}$.
Then $\varnothing \in \mathcal{S}$, since $\varnothing \cup B = B$, and if
$A \in \mathcal{S}$ and $s \in S$ then $A \cup \{s\} \in \mathcal{S}$,
since $A \cup B \in \mathcal{F}$ and so $(A \cup \{s\}) \cup B = (A \cup B) \cup \{s\} \in \mathcal{F}$.
Thus $\mathcal{S}$ is an inductive system, and in particular $\mathcal{F} \subset \mathcal{S}$.
Hence $A \cup B \in \mathcal{F}$ for each $A \in \mathcal{F}$.
\eop

This completes the proof of Theorem~\ref{theorem_free_11}. \eop

It should be clear that, if we had defined what it means for a set to be finite then
the least inductive system $\mathcal{F}$ is really just the set of finite subsets of $S$.
This fact is, however, irrelevant in the above proof.
Note that the usual definition of a set being finite requires properties of
the natural numbers which depend on all the Peano axioms, and so it can hardly be used in these notes.

%%% Local Variables: 
%%% mode: latex
%%% TeX-master: "sums"
%%% End: 

\startsection{Initial counting systems}

\label{initial}

In what follows the index set $S$ is considered to be fixed.
Counting systems are either $S$-typed counting systems or
counting systems with a single mapping, and the latter will
be referred to as single mapping counting systems.

We here present a characterisation of initial counting systems.
Theorem~\ref{theorem_initial_11} deals with the general case and is really just a translation
of the corresponding result for free commutative monoids.
However, the characterisation involves certain single mapping counting systems being initial
and Theorem~\ref{theorem_initial_11} itself gives no information about such counting systems.
This special case has to be dealt with separately, which is done in Theorems 
\ref{theorem_initial_21} and \ref{theorem_initial_31}. 
These results state that a
single mapping counting system is initial if and only if it is a Dedekind system, i.e., they consist
of the recursion theorem and its converse.

Let $(X,\Fam{f},x_0)$ be a counting system. For each $s \in S$ we define a new counting system
by replacing each $f_t$ with $t \ne s$ by $\id_X$. 
This results in the counting system $(X,\Fam{f}_s,x_0)$, where
$\Fam{f}_s = \{f_{s,t}\}_{t \in S}$  is the family of commuting mappings with
$f_{s,s} = f_s$ and $f_{s,t} = \id_X$ for all $t \ne s$. 
The counting system $(X,\Fam{f}_s,x_0)$ is really just a padded out version of the
single mapping counting system $(X,f_s,x_0)$, and we will also need the
minimal core $(X'_s,f'_s,x_0)$ of the $(X,f_s,x_0)$. Thus
$X'_s$ is the least $f_s$-invariant subset of $X$ containing $x_0$ and 
$f'_s$ is the restriction of $f_s$ to $X'_s$, considered as an element of $\Self{X'_s}$.

It should be clear that the minimal core of $(X,\Fam{f}_s,x_0)$ is just the counting system
$(X'_s,\Fam{f}'_s,x_0)$, where
$\Fam{f}'_s$ is obtained by restricting each of the mappings in $\Fam{f}_s$ to the set $X'_s$,
i.e.,
$\Fam{f}'_s = \{f'_{s,t}\}_{t \in S}$ with $f'_{s,s} = f'_s$ and $f'_{s,t} = \id_{X'_s}$ for all $t \ne s$.

Recall from Lemma~\ref{lemma_monoids_31} that
there is at most one morphism from a minimal counting system to any other counting system.

\begin{theorem}\label{theorem_initial_11}
Suppose $(X,\Fam{f},x_0)$ is minimal. Then the following are equivalent:

(1)\enskip 
The counting system $(X,\Fam{f},x_0)$ is initial.

(2)\enskip
For each $s \in S$ there exists a morphism $\pi_s : (X,\Fam{f},x_0) \to (X,\Fam{f}_s,x_0)$ and
the single mapping counting system $(X'_s,f'_s,x_0)$ is initial.

\end{theorem}

\proof 
(1) $\Rightarrow$ (2):
Fix $s \in S$; since $(X,\Fam{f},x_0)$ is initial there exists a morphism 
$\pi_s : (X,\Fam{f},x_0) \to (X,\Fam{f}_s,x_0)$.
Let $(Y,g,y_0)$ be a single mapping counting system and define a family
$\Fam{g} = \{g_s\}_{s \in S}$ of commuting mappings by letting $g_s = g$ and $g_t = \id_Y$ for all 
$t \ne s$. There thus exists a morphism $\pi : (X,\Fam{f},x_0) \to (Y,\Fam{g},y_0)$, since
$(X,\Fam{f},x_0)$ is initial, and in particular $\pi \circ f_s = g_s \circ \pi = g \circ \pi$. 
Let $\pi' : X'_s \to Y$ be the restriction of $\pi$ to $X'_s$. Then $\pi'(x_0) = y_0$ and
$\pi' \circ f'_s = g \circ \pi'$ and hence $\pi' : (X'_s,f'_s,x_0) \to (Y,g,y_0)$ is a morphism, and 
by Lemma~\ref{lemma_monoids_31} it is the unique morphism. This shows that
$(X'_s,f'_s,x_0)$ is initial.

(2) $\Rightarrow$ (1):
By Theorem~\ref{theorem_monoids_21} it is enough to 
show that the monoid associated with $(X,\Fam{f},x_0)$ is free with respect to the family $\{x_s\}_{s \in S}$.
For each $s \in S$ let $X_s$ be the least submonoid of $X$ containing $x_s$.

\begin{lemma}\label{lemma_initial_11}
$\ X_s = X'_s$ and
the submonoid $X_s$ is the monoid associated with the single mapping counting system
$(X'_s,f'_s,x_0)$. 
\end{lemma}

\proof 
Exactly as in the proof of Lemma~\ref{lemma_bam_31} $X_s$ is an $f_s$-invariant subset of
$X$ containing $x_0$ and so $X'_s \subset X_s$. Moreover,
as in the proof of Lemma~\ref{lemma_monoids_11}~(1) 
$X' = \{ x \in X : \mbox{$x + x' \in X'_s$ for all $x' \in X'_s$} \}$
is a submonoid of $X$ containing $x_s$ and thus $X_s \subset X'$. But $X' \subset X'_s$, since
$x_0 \in X'_s$, i.e., $X_s \subset X'_s$. This shows that $X_s = X'_s$.
Now $f'_s(x) = x_s + x$ for all $x \in X_s$, with $+$ the operation associated with $(X,\Fam{f},x_0)$
(and thus the operation on the submonoid $X_s$). Therefore
by the uniqueness in Theorem~\ref{theorem_mcs_11} 
$X_s$ is the monoid associated with the single mapping counting system
$(X'_s,f'_s,x_0)$. \eop

By Theorem~\ref{theorem_monoids_21} and Lemma~\ref{lemma_initial_11}
$X_s$ is free with respect to $x_s$ for each $s \in S$.

\begin{lemma}\label{lemma_initial_21}
If $(X,\Fam{f},x_0)$ is minimal and $\pi : (X,\Fam{f},x_0) \to (Y,\Fam{g},y_0)$
is a morphism then $\pi(X)$
is the least $\Fam{g}$-invariant subset $Y'$ of $Y$ containing $y_0$
and so $\pi$ can be considered as a morphism from 
$(X,\Fam{f},x_0)$ to the minimal core $(Y',\Fam{g}',y_0)$ of $(Y,\Fam{g},y_0)$.
\end{lemma}

\proof 
If $y = \pi(x) \in \pi(X)$ then
$g_s(y) = g_s(\pi(x)) = \pi(f_s(x))$ and so $g_s(y) \in \pi(X)$ for all $s \in S$. Thus $\pi(X)$ is a
$\Fam{g}$-invariant subset of $Y$ containing $y_0 = \pi(x_0)$, and hence
$Y' \subset \pi(X)$. On the other hand, $\pi^{-1}(Y')$
is an $\Fam{f}$-invariant subset of $X$ containing $x_0$ (since if
$\pi(x) \in Y'$ then $\pi(f_s(x)) = g_s(\pi(x)) \in Y'$) and so $\pi^{-1}(Y') = X$.
This implies that $\pi(X) \subset Y'$. 
\eop

By Lemma~\ref{lemma_initial_21} the morphism $\pi_s : (X,\Fam{f},x_0) \to (X,\Fam{f}_s,x_0)$ can be regarded
as a morphism $\pi_s : (X,\Fam{f},x_0) \to (X'_s,\Fam{f}'_s,x_0)$ and thus by Theorem~\ref{theorem_monoids_11}
$\pi_s : X \to X'_s$ is a homomorphism of the associated monoids
with $\pi_s(x_s) = x_s$ and $\pi_s(x_t) = x_0$ for all $t \ne s$.
As in the proof of Lemma~\ref{lemma_initial_11}
the monoid associated with $(X'_s,\Fam{f}'_s,x_0)$ is again the submonoid $X_s$ of $X$, which means that
$\pi_s$ can be considered as an endomorphism of $X$.
Therefore by Proposition~\ref{prop_free_31} the monoid $X$ is the internal direct sum of the
family $\{X_s\}_{s \in S}$ and, since
$X_s$ is free with respect to $x_s$ for each $s \in S$, Theorem~\ref{theorem_free_21} implies that
the monoid $X$ is free with respect to the family $\{x_s\}_{s \in S}$.

This completes the proof of Theorem~\ref{theorem_initial_11}. \eop

There are  two reasons why Theorem~\ref{theorem_initial_11} by itself
is not really satisfactory. The first is that it involves the single mapping counting systems 
$(X'_s,f'_s,x_0)$, $s \in S$, being initial and the second is that it gives no information about 
when this is the case. However, the characterisation of initial single mapping counting systems
is provided by the recursion theorem and its converse:
A single mapping counting system is initial if and only if it is a Dedekind system.
These results are given below.

In what follows a counting system now  always means a single mapping counting system.
Here is the recursion theorem:

\begin{theorem}\label{theorem_initial_21}
Each Dedekind system is initial.
\end{theorem}

\proof
In the proof we will need Lemma~\ref{lemma_spnat_21} several times, which we recall states that
if $(X,f,x_0)$ is minimal then for each $x \in X \setminus \{x_0\}$ there exists an $x' \in X$ such that
$x = f(x')$. This can also be applied to the minimal core; thus if
$(X,f,x_0)$ is any counting system and $X'$ is the least $f$-invariant subset of $X$ containing $x_0$
then for each $x \in X' \setminus \{x_0\}$ there exists an $x' \in X'$ such that $x = f(x')$.

Let $(X,f,x_0)$ be a Dedekind system and $(Y,g,y_0)$ be any counting system, and
consider the counting system $(X \times Y, f \times g, (x_0,y_0))$, where
$f \times g \in \Self{X \times Y}$ is given by
$(f \times g)(x,y) = (f(x),g(y))$ for all $x \in X$, $y \in Y$. Let $Z$ be the least
$(f \times g)$-invariant subset of $X \times Y$ containing $(x_0,y_0)$ and let
\[ X_0 = \{\, x \in X : \mbox{there exists exactly one $y \in Y$ such that $(x,y) \in Z$} \,\}\;. \]
We show that $X_0$ is an $f$-invariant subset of $X$ containing $x_0$, which implies that $X_0 = X$,
since by definition the Dedekind system $(X,f,x_0)$ is minimal.

The element $x_0$ is in $X_0$: Clearly
$(x_0,y_0) \in Z$, so suppose also $(x_0,y) \in Z$ for some $y \ne y_0$. Then 
$(x_0,y) \in Z \setminus \{(x_0,y_0)\}$ and hence by Lemma~\ref{lemma_spnat_21} there exists
$(x',y') \in Z$ with $(f \times g)(x',y') = (x_0,y)$.
In particular $f(x') = x_0$, which is not possible for the Dedekind system $(X,f,x_0)$. 
This shows that $x_0 \in X_0$.

Next let $x \in X_0$ and let $y$ be the unique element of $Y$ with $(x,y) \in Z$.
Hence $(f(x),g(y)) = (f \times g)(x,y) \in Z$, since $Z$ is $(f \times g)$-invariant.
Suppose also that $(f(x),y') \in Z$ for some $y' \in Y$.
Then $(f(x),y') \in Z \setminus \{(x_0,y_0)\}$, since $f(x) \ne x_0$, and 
so there exists $(x'',y'') \in Z$ with $(f(x''),g(y'')) = (f \times g)(x'',y'') = (f(x),y')$,
and in particular $f(x'') = f(x)$. But $f$ is injective, which means that $x'' = x$.
Thus $y'' = y$, since $x \in X_0$, which implies $y' = g(y'') = g(y)$.
This shows that $g(y)$ is the unique element of $Y$ with
$(f(x),g(y)) \in Z$ and in particular that $f(x) \in X_0$.

We have established that $X_0$ is an $f$-invariant subset of $X$ containing $x_0$, and so $X_0 = X$.
Now define a mapping $\pi : X \to Y$ by letting $\pi(x)$ be the unique element
of $Y$ such that $(x,\pi(x)) \in Z$ for each $x \in X$. 
Then $\pi(x_0) = y_0$, since $(x_0,y_0) \in Z$ and $\pi(f(x)) = g(\pi(x))$ for all $x \in X$,
since $(f(x),g(y)) \in Z$ whenever $(x,y) \in Z$ and so in particular
$(f(x),g(\pi(x))) \in Z$ for all $x \in X$. This gives us a morphism
$\pi : (X,f,x_0) \to (Y,g,y_0)$ which by
Lemma~\ref{lemma_monoids_31} is unique, and thus shows that $(X,f,x_0)$ is initial.
\eop

Here is the converse of the recursion theorem (Lawvere \cite{lawvere}):

\begin{theorem}\label{theorem_initial_31}
Each initial counting system is a Dedekind system.
\end{theorem}

\proof
We start with a simple general construction.
Let $(X,f,x_0)$ be any counting system,
choose an element $\omega$ not contained in $X$,
put $X_\omega = X \cup \{\omega\}$ and define a mapping $f_\omega : X_\omega \to X_\omega$ 
by letting $f_\omega(x) = f(x)$ for $x \in X$ and $f_\omega(\omega) = x_0$; 
thus $(X_\omega,f_\omega,\omega)$ is a counting system (and note that $f_\omega(X_\omega) \subset X$).

\begin{lemma}\label{lemma_initial_31}
If $(X,f,x_0)$ is minimal then so is $(X_\omega,f_\omega,\omega)$.
\end{lemma}

\proof 
Let $Y'$ be an $f_\omega$-invariant subset of $X_\omega$ containing $\omega$ and let
$Y = Y' \setminus \{\omega\}$. Then $x_0 = f_\omega(\omega) \in Y$ and
$Y$ is an $f$-invariant subset of $X$. (If $x \in Y$ then $f(x) = f_\omega(x) \in Y'$ and so
$f(x) \in Y$, since $f(x) \in X$.) Therefore $Y = X$, since
$(X,f,x_0)$ is minimal, and thus $Y' = X_\omega$. Hence $(X_\omega,f_\omega,\omega)$ is minimal.  
\eop

Suppose now $(X,f,x_0)$ is initial; then by Lemma~\ref{lemma_monoids_51}
$(X,f,x_0)$ is minimal and so by Lemma~\ref{lemma_initial_31} $(X_\omega,f_\omega,\omega)$ is also minimal.
Let $\pi : (X,f,x_0) \to (X_\omega,f_\omega,\omega)$ be the unique 
morphism. Consider the set $X_0 = \{ x \in X : f_\omega(\pi(x)) = x \}$; then $x_0 \in X_0$, since
$f_\omega(\pi(x_0)) = f_\omega(\omega) = x_0$ and if $x \in X_0$ then $f_\omega(\pi(x)) = x$ and so
\[  f_\omega(\pi(f(x))) = f_\omega(f_\omega(\pi(x))) = f_\omega(x) = f(x)\;, \]
i.e., $f(x) \in X_0$. Thus $X_0$ is an $f$-invariant subset of $X$ containing $x_0$ and hence
$X_0 = X$, since $(X,f,x_0)$ is minimal. 
This means that
$f_\omega \circ \pi = \id_X$ with $f_\omega$ here considered as a mapping $X_\omega \to X$.
Therefore $f_\omega: X_\omega \to X_\omega$ is injective but by
definition $\omega \notin f_\omega(X_\omega)$; moreover, we have seen above that
$(X_\omega,f_\omega,\omega)$ is minimal, i.e., $(X_\omega,f_\omega,\omega)$ is a Dedekind system.

This implies that $(X,f,x_0)$ is a Dedekind system:
By Theorem~\ref{theorem_initial_21} $(X_\omega,f_\omega,\omega)$ is initial and so by Lemma~\ref{lemma_monoids_41} 
the two initial counting systems $(X,f,x_0)$ and $(X_\omega,f_\omega,\omega)$ are isomorphic. By
Lemma~\ref{lemma_initial_41}~(2) below $(X,f,x_0)$ is then a Dedekind system, since it is isomorphic to the
Dedekind system $(X_\omega,f_\omega,\omega)$. \eop

\begin{lemma}\label{lemma_initial_41}
Let $(X,f,x_0)$ and $(Y,g,y_0)$ be isomorphic counting systems.

(1)\enskip
If $(X,f,x_0)$ is minimal then so is $(Y,g,y_0)$.

(2)\enskip
If $(X,f,x_0)$ is a Dedekind system then so is $(Y,g,y_0)$.
\end{lemma}

\proof 
Let $\pi : (X,f,x_0) \to (Y,g,y_0)$ be an isomorphism. 

(1)\enskip
Let $Y' \subset Y$ be a $g$-invariant subset of $Y$ containing $y_0$ and put $X' = \pi^{-1}(Y')$.
Then $x_0 \in X'$ (since $\pi(x_0) = y_0 \in Y_0$) and 
$X'$ is $f$-invariant (since if $x \in X'$ then $\pi(x) \in Y'$, hence
$\pi(f(x)) = g(\pi(x)) = g(y) \in Y'$ and so $f(x) \in X'$). Therefore
$X' = X$, which implies that
$Y' = \pi(\pi^{-1}(Y')) = \pi(X') = \pi(X) = Y$. This shows that $(Y,g,y_0)$ is minimal.

(2)\enskip
By (1) $(Y,g,y_0)$ is minimal.
Moreover, the mapping $g = \pi \circ f \circ \pi^{-1}$, as the composition of three injective mappings, 
is itself injective.
Finally, if $y_0 = g(y)$ for some $y \in Y$ then 
$\pi(x_0) = y_0 = g(y) = g(\pi(\pi^{-1}(y)) = \pi(f(\pi^{-1}(y)))$
and thus  $x_0 = f(\pi^{-1}(y)) \in f(X)$. But this is not the case and hence 
$y_0 \notin g(Y)$. \eop

The proof of Theorem~\ref{theorem_initial_31} first shows that the existence of an initial
counting system implies that of a Dedekind system, and then uses the Dedekind system
to show that any initial counting system is a Dedekind system.
If one is prepared to accept that a Dedekind system exists then the first step is not needed and
the proof is then shorter.
It is easy to see that the existence of a Dedekind system is equivalent to that of a Dedekind-infinite set,
i.e., a set $X$ for which there exists an injective mapping $f : X \to X$ which is not surjective.
(If $f : X \to X$ is such a mapping and $x_0 \notin f(X)$ then the minimal core of the
counting system $(X,f,x_0)$ is a Dedekind system.) Thus in a world where all sets are finite
there are neither initial counting systems nor Dedekind systems.

\bigskip
\bigskip

{\sc Fakult\"at f\"ur Mathematik, Universit\"at Bielefeld}\\
{\sc Postfach 100131, 33501 Bielefeld, Germany}\\
\textit{E-mail address:} \texttt{preston@math.uni-bielefeld.de}\\
\textit{URL:} \texttt{http://www.math.uni-bielefeld.de/\symbol{126}preston}

%%% Local Variables: 
%%% mode: latex
%%% TeX-master: "sums"
%%% End: 

\end{document}